\documentclass[11pt]{article}
\usepackage{amsmath, amssymb, amscd, epsfig, amsthm}

\newtheorem{thm}{Theorem}
\newtheorem{prop}[thm]{Proposition}
\newtheorem{lem}[thm]{Lemma}

\newtheorem{cor}[thm]{Corollary}
\newtheorem{rem}[thm]{Remark}
\newtheorem{df}[thm]{Definition}
\newtheorem{ex}[thm]{Example}

\renewcommand{\epsilon}{\varepsilon}
\renewcommand{\phi}{\varphi}
\renewcommand{\deg}{\operatorname{deg}}

\newcommand{\BB}{\mathbb}

\newcommand{\separate}{\vskip5pt}

\newcommand{\re}{\operatorname{Re}}
\newcommand{\tr}{\operatorname{Tr}}
\newcommand{\End}{\operatorname{End}}

\newcommand{\HC}{\BB H_{\BB C}}

\newcommand{\C}{\operatorname{C}}

\usepackage[OT2,T1]{fontenc}

\begin{document}

\title{\bf Schwarz Lemma in Quaternionic Setting and $\BB C^4$}
\author{Chris Connell and Matvei Libine}
\maketitle

\begin{abstract}
We study a Schwarz Lemma on the matrix ball
\[
\BB D^+=\{Z\in M_2(\BB C);\, \|Z\|<1\}
       \simeq U(2,2)/(U(2)\times U(2)),
\]
which arises naturally from the complexification of the quaternions.  We
prove Schwarz and Schwarz--Pick estimates on $\BB D^+$, construct its
$U(2,2)$-invariant Riemannian Poincar\'e metric, and describe a Cayley
transform to a tube-domain realization.  The intrinsic spectral metric is
the Kobayashi metric and is contracted by holomorphic maps.  Comparing it
with the Poincar\'e metric gives the sharp estimate $d(F(P),F(Q))\le \sqrt2\,d(P,Q)$ for every holomorphic self-map $F$ of $\BB D^+$; in general the constant
$1$ is false.  We also prove an infinitesimal rigidity statement and record
a Burns--Krantz boundary rigidity theorem at the smooth part of
$\partial\BB D^+$.
Finally, we develop a conformally natural barycenter construction on the
Shilov boundary $U(2)$.  Every admissible probability measure has a unique
barycenter in $\BB D^+$, the construction is $U(2,2)$-natural, and it gives
a barycentric extension of suitable boundary maps.  Under the natural
admissibility hypothesis, a holomorphic self-map extending to $U(2)$ is the
barycentric extension of its boundary values.
Many aspects of this paper are motivated by quaternionic analysis,
however, no prior knowledge of quaternionic analysis is required from
the reader.
\end{abstract}


\section{Introduction}

Let ${\cal D}=\{z\in\BB C:|z|<1\}$ be the unit disk. The classical
Schwarz Lemma asserts that a holomorphic map $f:{\cal D}\to{\cal D}$ with
$f(0)=0$ satisfies $|f(z)|\le |z|$. Its invariant form is the
Schwarz--Pick theorem, and its geometric extensions include the
Schwarz--Pick--Ahlfors and Yau Schwarz lemmas; see \cite{A,O,Kr2,Y}.
Some generalizations of the Schwarz lemma in the context of Clifford analysis appear in \cite{YQ,Z1,Z2};
these results readily translate to quaternionic analysis.

A direct replacement of holomorphic functions by quaternionic regular
functions on the Euclidean unit ball in $\BB H$ does not satisfy the
expected Schwarz estimate, as the counterexample of Section 3 shows. This
leads us instead to the domain
\[
\BB D^+=\{Z\in M_2(\BB C):\|Z\|<1\}
             \simeq U(2,2)/(U(2)\times U(2)),
\]
where $\|\cdot\|$ is the spectral norm\footnote{In contrast, the unit disk ${\cal D} = U(1,1)/(U(1) \times U(1))$.}.
Although this is the classical
matrix ball of type $I_{2,2}$, it arises naturally from quaternionic
analysis: the complexification $\BB H\otimes_{\BB R}\BB C$ is
$\HC \cong M_2(\BB C)$, and regular functions on the quaternionic unit ball extend
holomorphically to $\BB D^+$.

We prove Schwarz and Schwarz--Pick estimates for holomorphic maps on
$\BB D^+$, describe its fractional-linear $U(2,2)$ action, construct its
invariant Riemannian Poincar\'e metric, and give a Cayley transform to a
tube-domain realization. The Riemannian metric is built from the
Hilbert--Schmidt norm. Replacing that norm by the spectral norm gives an
invariant Finsler metric, which we identify with the Kobayashi (and
Carath\'eodory) metric. Every holomorphic self-map contracts this intrinsic
metric, and comparison of the two metrics yields the sharp estimate
\[
 d(F(P),F(Q))\le\sqrt2\,d(P,Q).
\]
The factor $\sqrt2$ is attained, so the Riemannian Poincar\'e metric is not
in general distance decreasing.  The same point of view also gives a useful
rigidity statement: if $F(0)=0$ and $dF_0$ preserves the infinitesimal
Kobayashi norm in every tangent direction, then $F$ is one of the standard
linear automorphisms $Z\mapsto UZV$ or $Z\mapsto UZ^tV$.  At a smooth point
of the Euclidean boundary we also obtain a boundary Schwarz lemma directly
from the Burns--Krantz rigidity theorem of Ng--Rong \cite{NR}.

We also develop a conformally natural barycenter construction on the Shilov
boundary $U(2)$, in the spirit of Douady--Earle \cite{DE} and
Besson--Courtois--Gallot \cite{BCG}. An admissible probability measure, 
namely one assigning mass less than $1/2$ to every nontransversality set
${\cal Z}_W$, has a unique barycenter in $\BB D^+$. The assignment is
$U(2,2)$-natural and gives a barycentric extension of suitable boundary
maps. Under the admissibility hypothesis of Theorem~\ref{F-is-extension-thm},
a holomorphic self-map extending continuously to the Shilov boundary is the
barycentric extension of its own boundary values. In particular, the
normalization $F(0)=0$ can be detected from a boundary moment condition.

We want to emphasize that all results of this paper -- including our Schwarz and Schwarz--Pick estimates and the barycenter construction on $\BB D^+$ -- are completely natural from the point of view of the conformal group action and that the open set $\BB D^+$ is distinct from everpresent Euclidean balls or polydisks.

Sections 2--7 establish the quaternionic and symmetric-space background and
the basic Schwarz results. Section~\ref{spectral-sect} proves the sharp
metric comparison and the interior and smooth-boundary rigidity statements.
Section~\ref{bary-sect} develops the barycenter construction and its
applications.  We also point out there that the higher-rank barycenter
method of Connell--Farb applies to compact locally symmetric quotients of
$\BB D^+$, giving the corresponding degree and minimal-volume consequences.
We conclude the paper with a list of further
questions arising naturally from this version of the Schwarz Lemma (Section \ref{questions-sect}).

\subsection*{Acknowledgements.}
The first author gratefully acknowledges partial support from the Simons Foundation (grant \#965245) and the National Science Foundation (grant DMS-2514510).

\paragraph{\bf Use of artificial intelligence.}
OpenAI's ChatGPT-5 was used for light editorial assistance, and for most of the applications in Sections \ref{spectral-sect} and \ref{bary-sect} after 8.2, but not the overall theme of those applications within each of those two sections. Everything before Section \ref{spectral-sect} is original to the authors. The authors independently verified all arguments and computations and take full responsibility for the content.

\section{Preliminaries}\label{preliminaries}

In this section we establish notations and state relevant results from
quaternionic analysis.
A contemporary review of quaternionic analysis can be found in \cite{CSSS, Su}.
Quaternionic analysis also has many applications in physics
(see, for instance, \cite{GT}).

\subsection{Complexified Quaternions $\HC$ and the Conformal Group $GL(2,\HC)$}

We use notations established in \cite{FL1}.
In particular, ${\bf e_0}$, ${\bf e_1}$, ${\bf e_2}$, ${\bf e_3}$
denote the units of the classical quaternions $\BB H$ corresponding to the
more familiar $1$, $i$, $j$, $k$
(we reserve the symbol $i$ for $\sqrt{-1} \in \BB C$).
Thus, $\BB H$ is an algebra over $\BB R$ generated by
${\bf e_0}$, ${\bf e_1}$, ${\bf e_2}$, ${\bf e_3}$,
and the multiplicative structure is determined by the rules
\[
{\bf e_0 e_i} = {\bf e_i e_0} = {\bf e_i}, \quad
({\bf e_i})^2 = {\bf e_1 e_2 e_3} = - {\bf e_0}, \quad
{\bf e_ie_j} = - {\bf e_je_i}, \quad 1 \le i< j \le 3,
\]
and the fact that $\BB H$ is a division ring.
In this paper we mostly consider the algebra of complexified quaternions
$\HC = \BB C \otimes \BB H$ (also known as biquaternions),
it can be identified with the algebra of $2 \times 2$ complex matrices:
\begin{multline*}
\HC = \BB H \otimes_{\BB R} \BB C
= \{ Z = z_0{\bf e_0} + z_1{\bf e_1} + z_2{\bf e_2} + z_3{\bf e_3} ;\:
z_0,z_1,z_2,z_3 \in \BB C \}  \\
\simeq \biggl\{ Z= \begin{pmatrix} z_0-iz_3 & -iz_1-z_2 \\ -iz_1+z_2 & z_0+iz_3
\end{pmatrix} ; \: z_k \in \BB C \biggr\}
= \biggl\{ Z = \begin{pmatrix} z_{11} & z_{12} \\ z_{21} & z_{22} \end{pmatrix}
; \: z_{ij} \in \BB C \biggr\}.
\end{multline*}
We denote by $\BB S$ (respectively $\BB S'$) the irreducible 2-dimensional
left (respectively right) $\HC$-module, as described in Subsection 2.3
of \cite{FL1}. The spaces $\BB S$ and $\BB S'$ can be realized as respectively
columns and rows of complex numbers:
\[
\BB S = \bigl\{ \bigl(\begin{smallmatrix} z_1 \\ z_2 \end{smallmatrix}\bigr);\:
z_1, z_2 \in \BB C \bigr\}, \qquad
\BB S' = \{ (z_1, z_2);\: z_1, z_2 \in \BB C \}.
\]
For
$Z = z_0{\bf e_0} + z_1{\bf e_1} + z_2{\bf e_2} + z_3 {\bf e_3} \in \HC$,
we denote by $Z^+$ its {\em quaternionic conjugate}:
\[
Z^+ = z_0{\bf e_0} - z_1{\bf e_1} - z_2{\bf e_2} - z_3 {\bf e_3}
\]
(the bar notation is reserved for complex conjugation), and write
\begin{multline*}
N(Z) = ZZ^+ = Z^+Z  =
\det \begin{pmatrix} z_{11} & z_{12} \\ z_{21} & z_{22} \end{pmatrix}
= z_{11}z_{22}-z_{12}z_{21} \\
= (z_0)^2 + (z_1)^2 + (z_2)^2 + (z_3)^2
\quad \in \BB C
\end{multline*}
and think of it as (the square of) the norm of $Z$.
For $X = x_0{\bf e_0} + x_1{\bf e_1} + x_2{\bf e_2} + x_3{\bf e_3} \in \BB H$,
let
\[
\|X\|_{\BB H} = \sqrt{N(X)} = \sqrt{(x_0)^2+(x_1)^2+(x_2)^2+(x_3)^2}
\]
be the Euclidean norm of $X$.
An element $Z$ in $\HC$ (or $\BB H$) has a multiplicative inverse
$Z^{-1}$ if and only if $N(Z) \ne 0$. Let
\[
\HC^{\times} = \{ Z \in \HC ;\: N(Z) \ne 0 \}
= \{\text{invertible complexified quaternions}\}.
\]
We realize the group of the unitary matrices $U(2)$ as
\[
U(2) = \{ Z \in \HC ;\: Z^*=Z^{-1} \},
\]
where $Z^*$ denotes the complex conjugate transpose of a complex matrix
$Z$. For $R>0$, we set
\[
U(2)_R = \{ RZ ;\: Z \in U(2) \} \quad \subset \HC.
\]

We review the action of a group $GL(2,\HC) \simeq GL(4,\BB C)$ on $\HC$
by fractional linear (or conformal) transformations.
Firstly, the tautological action of $GL(2,\HC)$ on $\HC=\BB C^4$ extends
to an action on the complex Grassmannian $\operatorname{Gr}_2(\HC)$
consisting of complex $2$-dimensional planes inside $\HC$.
The elements of $\operatorname{Gr}_2(\HC)$ can be realized using two linearly
independent column vectors in $\BB C^4$ that we can combine into
$4 \times 2$ matrices
$\bigl(\begin{smallmatrix} X \\ Y \end{smallmatrix}\bigr)$, $X, Y \in \HC$,
with two columns $\bigl(\begin{smallmatrix} X \\ Y \end{smallmatrix}\bigr)$
and $\bigl(\begin{smallmatrix} X' \\ Y' \end{smallmatrix}\bigr)$
describing the same element of $\operatorname{Gr}_2(\HC)$ if and only if
they span the same plane:
\[
\begin{pmatrix} X \\ Y \end{pmatrix}
\sim \begin{pmatrix} X' \\ Y' \end{pmatrix}
\qquad \Longleftrightarrow \qquad
\begin{pmatrix} X' \\ Y' \end{pmatrix}
= \begin{pmatrix} Xa \\ Ya \end{pmatrix} \qquad
\text{for some $a \in \HC^{\times}$}.
\]
Then elements $Z \in \HC$ can be identified with planes in
$\operatorname{Gr}_2(\HC)$ that have a basis of the form
$\bigl(\begin{smallmatrix} Z \\ 1 \end{smallmatrix}\bigr)$.
This identifies $\HC$ with an open dense subset of $\operatorname{Gr}_2(\HC)$.
The action of $GL(2,\HC)$ on these elements is
\[
\begin{pmatrix} a & b \\ c & d \end{pmatrix}
\begin{pmatrix} Z \\ 1 \end{pmatrix}
= \begin{pmatrix} aZ+b \\ cZ+d \end{pmatrix}
\sim \begin{pmatrix} (aZ+b)(cZ+d)^{-1} \\ 1 \end{pmatrix},
\qquad \begin{pmatrix} a & b \\ c & d \end{pmatrix} \in GL(2,\HC).
\]
(Of course, such identification is possible only if $cZ+d$ is invertible.)
In practice, we restrict ourselves to $\HC$ and express this action of
$GL(2,\HC)$ using the fractional linear (or conformal) transformations that
may be undefined on a subset of measure 0:
\begin{equation}  \label{conformal-action}
h: Z \mapsto (aZ+b)(cZ+d)^{-1} = (a'-Zc')^{-1}(-b'+Zd'),
\qquad Z \in \HC,
\end{equation}
where
$h = \bigl(\begin{smallmatrix} a & b \\ c & d \end{smallmatrix}\bigr)
\in GL(2,\HC)$ and 
$h^{-1} = \bigl(\begin{smallmatrix} a' & b' \\ c' & d' \end{smallmatrix}\bigr)$.

\subsection{Harmonic and Regular Functions on $\BB H$}

The starting point of quaternionic analysis is the notion of a
(left or right) regular function.
Such functions are defined using a quaternionic analogue of the
Cauchy-Riemann equations. Write $X \in \BB H$ as
$X = x_0{\bf e_0} + x_1{\bf e_1} + x_2{\bf e_2} + x_3{\bf e_3}$,
$x_0, x_1, x_2, x_3 \in \BB R$,
and factor the four-dimensional Laplacian operator $\square$ on $\BB H$
as a product of two Dirac operators
\[
\square = \frac {\partial^2}{(\partial x_0)^2} +
\frac {\partial^2}{(\partial x_1)^2} +
\frac {\partial^2}{(\partial x_2)^2} +
\frac {\partial^2}{(\partial x_3)^2}
= \nabla \nabla^+ = \nabla^+ \nabla,
\]
where
\begin{align*}
\nabla^+ &= {\bf e_0} \frac{\partial}{\partial x_0}
+ {\bf e_1} \frac{\partial}{\partial x_1}
+ {\bf e_2} \frac{\partial}{\partial x_2}
+ {\bf e_3} \frac{\partial}{\partial x_3} \qquad \text{and} \\
\nabla &= {\bf e_0} \frac{\partial}{\partial x_0}
- {\bf e_1} \frac{\partial}{\partial x_1}
- {\bf e_2} \frac{\partial}{\partial x_2}
- {\bf e_3} \frac{\partial}{\partial x_3}.
\end{align*}
The operators $\nabla^+$, $\nabla$ can be applied to functions
on the left and on the right. Let $U \subset \BB H$ be an open subset.
A differentiable function $\phi$ on $U$ (usually with values in $\BB R$,
$\BB C$, $\BB H$, $\BB S$, $\BB S'$ or $\HC$) is called {\em harmonic}
if $(\square \phi)(X)=0$ for all $X \in U$.
A differentiable function $f$ on $U$ with values in $\BB H$, $\BB S$ or $\HC$
is called {\em left-regular} if
\[
(\nabla^+ f)(X) =
{\bf e_0} \frac{\partial f}{\partial x_0}(X)
+ {\bf e_1} \frac{\partial f}{\partial x_1}(X)
+ {\bf e_2} \frac{\partial f}{\partial x_2}(X)
+ {\bf e_3} \frac{\partial f}{\partial x_3}(X) =0
\]
for all $X \in U$.
Similarly, a differentiable function $g$ on $U$ with values in $\BB H$,
$\BB S'$ or $\HC$ is {\em right-regular} if
\[
(g \nabla^+)(X) =
\frac{\partial g}{\partial x_0}(X) {\bf e_0}
+ \frac{\partial g}{\partial x_1}(X) {\bf e_1}
+ \frac{\partial g}{\partial x_2}(X) {\bf e_2}
+ \frac{\partial g}{\partial x_3}(X) {\bf e_3} =0
\]
for all $X \in U$.

We restate the Poisson formula for harmonic functions as in
Theorem 34 in \cite{FL1}:

\begin{thm} \label{Poisson-thm}
Let $\phi$ be a real analytic harmonic function defined on a neighborhood of
the closure of a ball ${\cal B}_R = \{X \in \BB H ;\: \|X\|_{\BB H} < R \}$
of radius $R>0$. Then
\[
\phi(X_0) = \frac 1{2\pi^2} \int_{X \in \partial {\cal B}_R}
\frac{R^2 - N(X_0)}{N(X-X_0)^2} \cdot \frac{dS}R \cdot \phi(X),
\qquad \text{for all } X_0 \in {\cal B}_R.
\]
\end{thm}

As usual, $\partial {\cal B}_R$ denotes the boundary sphere of ${\cal B}_R$,
and $dS$ is the standard Euclidean volume element on $\partial {\cal B}_R$.

Similarly, the left and right regular functions satisfy an analogue
of Cauchy's integral formula called Cauchy-Fueter formulas.
We introduce a quaternionic valued 3-form on $\BB H$
\[
Dx= {\bf e_0} dx^1 \wedge dx^2 \wedge dx^3
- {\bf e_1} dx^0 \wedge dx^2 \wedge dx^3
+ {\bf e_2} dx^0 \wedge dx^1 \wedge dx^3
- {\bf e_3} dx^0 \wedge dx^1 \wedge dx^2.
\]

\begin{thm} [Cauchy-Fueter Formulas, \cite{F1,F2}]  \label{Fueter}
Let $U \subset \BB H$ be an open bounded subset with piecewise ${\cal C}^1$
boundary $\partial U$. Suppose that $f(X)$ is left-regular on a
neighborhood of the closure $\overline{U}$, then
\[
\frac1{2\pi^2} \int_{\partial U} \frac{(X-X_0)^{-1}}{N(X-X_0)} \cdot Dx \cdot f(X)
= \begin{cases}
f(X_0) & \text{if $X_0 \in U$;} \\
0 & \text{if $X_0 \notin \overline{U}$.}
\end{cases}
\]
If $g(X)$ is right-regular on a neighborhood of the closure $\overline{U}$,
then
\[
\frac1{2\pi^2} \int_{\partial U} g(X) \cdot Dx \cdot \frac {(X-X_0)^{-1}}{N(X-X_0)}
= \begin{cases}
g(X_0) & \text{if $X_0 \in U$;} \\
0 & \text{if $X_0 \notin \overline{U}$.}
\end{cases}
\]
\end{thm}

For the orientation of the contours of integration we refer to the discussion
before Lemma 6 in \cite{FL1}.
When the open set $U$ is the open ball
${\cal B}_R = \{X \in \BB H ;\: \|X\|_{\BB H} < R \}$
of radius $R>0$ centered at the origin, the boundary
$\partial U$ is the sphere
$\partial {\cal B}_R = \{X \in \BB H ;\: \|X\|_{\BB H} = R \}$,
and the restriction
\begin{equation}  \label{Dx-restriction}
Dx \Bigr|_{\partial {\cal B}_R} = \frac{X}R\,dS,
\end{equation}
where $dS$ is the Euclidean volume element on $\partial {\cal B}_R$ as before
(Lemma 6 in \cite{FL1}).

\subsection{Harmonic and Regular Functions on $\HC$}

Now let $U$ be an open subset of $\HC$.
As in Section 2 of \cite{split}, we consider functions on $U$ that are
holomorphic with respect to the complex variables $z_0$, $z_1$, $z_2$, $z_3$
(or, equivalently, $z_{11},z_{12},z_{21},z_{22}$).
A holomorphic function $\phi$ on $U$ (usually with values in $\BB C$,
$\BB S$, $\BB S'$ or $\HC$) is called {\em holomorphic harmonic}
if $(\square \phi)(Z)=0$ for all $Z \in U$.
This is equivalent to saying that $\phi$ is holomorphic and annihilated by
a differential operator
\[
4(\partial_{11}\partial_{22} - \partial_{12}\partial_{21})
= \frac{\partial^2}{(\partial z^0)^2} + \frac{\partial^2}{(\partial z^1)^2}
+ \frac{\partial^2}{(\partial z^2)^2} + \frac{\partial^2}{(\partial z^3)^2},
\qquad \text{where } \partial_{ij}= \frac{\partial}{\partial z_{ij}}.
\]
A holomorphic function $f$ on $U$ with values in $\BB S$ or $\HC$
is called {\em holomorphic left-regular} if $(\nabla^+ f)(Z) =0$
for all $Z \in U$.
Similarly, a holomorphic function $g$ on $U$ with values in $\BB S'$ or $\HC$
is {\em holomorphic right-regular} if $(g \nabla^+)(Z) =0$ for all $Z \in U$.

For holomorphic functions $f$ with values in $\BB S$ or $\HC$ and
$g$ with values in $\BB S'$ or $\HC$, the following derivatives are equal:
\[
\nabla^+ f =
{\bf e_0} \frac{\partial f}{\partial z^0} 
+ {\bf e_1} \frac{\partial f}{\partial z^1}
+ {\bf e_2} \frac{\partial f}{\partial z^2}
+ {\bf e_3} \frac{\partial f}{\partial z^3}
= 2 \begin{pmatrix} \partial_{22} & - \partial_{21}  \\
- \partial_{12} & \partial_{11} \end{pmatrix} f,
\]
\[
g \nabla^+ =
\frac{\partial g}{\partial z^0} {\bf e_0}
+ \frac{\partial g}{\partial z^1} {\bf e_1}
+ \frac{\partial g}{\partial z^2} {\bf e_2}
+ \frac{\partial g}{\partial z^3} {\bf e_3}
=
2g \begin{pmatrix} \partial_{22} & - \partial z_{21}  \\
- \partial z_{12} & \partial_{11} \end{pmatrix}.
\]

\section{A Counterexample}

Let $\cal D$ be the open unit disk in $\BB C$ and $\cal B$ the open unit ball
in $\BB H$:
\[
{\cal D} = \{ z \in \BB C ;\: |z| <1 \}, \qquad
{\cal B} = \{ X \in \BB H ;\: \|X\|_{\BB H} <1 \};
\]
then we let $\overline{\cal D}$ and $\overline{\cal B}$ denote their respective
closures.
In this section we replace a holomorphic function $f$ on $\cal D$ with
a left regular function $F$ on $\cal B$ and construct a counterexample to
the most obvious statement resembling Schwarz Lemma.
Consider a holomorphic function
\[
f(z)=3z-z^3 = (3x-x^3+3xy^2) + i(3y+y^3-3x^2y), \qquad z=x+iy,
\]
and ``replicate'' it as a left regular function in three ways:
\[
F_i({\bf e_0}x_0 + {\bf e_1}x_1 + {\bf e_2}x_2 + {\bf e_3}x_3) =
{\bf e_0} \bigl( 3x_0-(x_0)^3+3x_0(x_i)^2 \bigr)
+ {\bf e_i} \bigl( 3x_i+(x_i)^3-3(x_0)^2x_i \bigr),
\]
$i=1,2,3$. Set
\begin{multline*}
F(X)=F_1(X)+F_2(X)+F_3(X)
= {\bf e_0} 3x_0 \bigl( 3-(x_0)^2+ (x_1)^2 + (x_2)^2 + (x_3)^2 \bigr)  \\
+ {\bf e_1} x_1 \bigl( 3+(x_1)^2-3(x_0)^2 \bigr)
+ {\bf e_2} x_2 \bigl( 3+(x_2)^2-3(x_0)^2 \bigr)
+ {\bf e_3} x_3 \bigl( 3+(x_3)^2-3(x_0)^2 \bigr).
\end{multline*}
Clearly, $F_1$, $F_2$, $F_3$, $F$ are left regular functions $\BB H \to \BB H$,
\[
F(0)=0 \qquad \text{and} \qquad
\frac{\partial F}{\partial x_0} \biggr|_{X=0} = 9.
\]
Thus, the quotient
\[
\frac{\|F(X)\|_{\BB H}}{\|X\|_{\BB H}} \quad
\text{can be made arbitrary close to $9$}
\]
when $X$ is near $0$.
On the other hand, we can find the maximum of $\|F(X)\|_{\BB H}$ over
$X \in \overline{\cal B}$ and see that it is strictly less than $9$.
While it is not hard to find the exact maximum value of $\|F(X)\|_{\BB H}$,
it is easier to estimate it from above.

\begin{lem}
We have: 
\[
\sup_{X \in \overline{\cal B}} \|F(X)\|_{\BB H} \le 4\sqrt{3} < 7.
\]
\end{lem}

\begin{proof}
We observe that the absolute value of the ${\bf e_0}$-component of $F$ is
\[  
3|x_0| \bigl( 3-(x_0)^2+ (x_1)^2 + (x_2)^2 + (x_3)^2 \bigr)
\le 3|x_0| \bigl( 4 - 2(x_0)^2 \bigr),
\]
since $(x_1)^2 + (x_2)^2 + (x_3)^2 \le 1 - (x_0)^2$.
On the other hand, for each $i=1,2,3$, the absolute value of the
${\bf e_i}$-component of $F$ is
\[
|x_i| \bigl( 3+(x_i)^2-3(x_0)^2 \bigr)
\le 4 \sqrt{1-(x_0)^2} \bigl( 1 - (x_0)^2 \bigr),
\]
since $(x_i)^2 \le 1 - (x_0)^2$.
Hence we can estimate
\[
\|F(X)\|_{\BB H}^2 \le 36 (x_0)^2 \bigl( 2-(x_0)^2 \bigr)^2
+ 48 \bigl( 1-(x_0)^2 \bigr)^3
= 12 \bigl( 4 - (x_0)^6 \bigr)
\le 48,
\]
for all $X \in \overline{\cal B}$.
\end{proof}  

Rescaling $F(X)$ we conclude:

\begin{rem}
There exists a left regular function $F: \BB H \to \BB H$ such that
\[
F(0)=0, \qquad \|F(X)\|_{\BB H} \le 1 \text{ for all $X \in \overline{\cal B}$}
\qquad \text{and} \qquad \frac{\|F(X_0)\|_{\BB H}}{\|X_0\|_{\BB H}} \ge 1
\]
for some $X_0 \in {\cal B}$, $X_0 \ne 0$.

A similar observation can be made for right regular functions
$G: \BB H \to \BB H$.
\end{rem}

\section{Quaternionic Analogue of the Unit Disk}

One of the main points of this paper is that the quaternionic analogue of
the unit disk ${\cal D} = \{ z \in \BB C ;\: |z| <1 \}$ is not the Euclidean
unit ball ${\cal B} = \{ X \in \BB H ;\: \|X\|_{\BB H} <1 \}$ but rather
a certain open domain $\BB D^+ \subset \HC$.

\subsection{Open Domains $\BB D^+$ and $\BB D^-$}

We consider the complexified quaternions $\HC = \BB H \otimes_{\BB R} \BB C$
and treat them as $2 \times 2$ complex matrices.
Recall a special case of the singular value decomposition for such matrices:

\begin{prop}  \label{svd}
Every element $Z \in \HC$ can be written as a product $Z=UDV^*$ with
$U, V \in U(2)$ and
\begin{equation}  \label{Diag}
  D = \bigl(\begin{smallmatrix} \lambda_1 & 0 \\
    0 & \lambda_2 \end{smallmatrix}\bigr) \in \HC,
\qquad \lambda_1, \lambda_2 \in \BB R, \: \lambda_1, \lambda_2 \ge 0.
\end{equation}
\end{prop}

For a $2 \times 2$ complex matrix $A \in \HC$, we define the
{\em (spectral) norm}
\begin{equation}  \label{spectral-norm}
\|A\| = \sup\{ \|Az\|;\: z \in \BB C^2,\: \|z\|=1\}.
\end{equation}
(If $z= \bigl(\begin{smallmatrix} z_1 \\ z_2 \end{smallmatrix}\bigr)
\in \BB C^2$, let $\|z\|=\sqrt{|z_1|^2+|z_2|^2}$.)
Clearly, if $A$ is written as a product $A=UDV^*$ as in Proposition \ref{svd},
then
\[
\|A\| = \max\{\lambda_1, \lambda_2\}.
\]

\begin{lem}
If $A \in \BB H \subset \HC$, then the spectral norm
\begin{equation}  \label{equal-norms-eqn}
\|A\| = \sqrt{N(A)} = \|A\|_{\BB H}.
\end{equation}
\end{lem}  

\begin{proof}
If $A \in \BB H$, then either $A=0$ or $A$ can be written as
$A = \sqrt{N(A)} \cdot U$ with $U = N(A)^{-1/2}A \in SU(2)$.
In either case we obtain \eqref{equal-norms-eqn}.
\end{proof}

We address the uniqueness of the matrix $D$ in the decomposition
\eqref{Diag}.
For any $Z \in \HC$, the matrix $ZZ^*$ is a $2 \times 2$ self-adjoint
complex matrix.
If $Z$ is written as a product $Z=UDV^*$ as in Proposition \ref{svd},
$ZZ^* = UD^2U^*$. Hence $ZZ^*$ is conjugate by an element of $U(2)$ to
$D^2 = \bigl(\begin{smallmatrix} \lambda_1^2 & 0 \\
  0 & \lambda_2^2 \end{smallmatrix}\bigr)$,
and $\|ZZ^*\|=\|Z\|^2$.
Since the diagonal entries $\lambda_1$ and $\lambda_2$ in \eqref{Diag}
are the square roots of the eigenvalues of $ZZ^*$, the pair
$\{\lambda_1,\lambda_2\}$ is uniquely determined by $Z$.
(Matrix $D$, however, typically is {\em not} unique because the entries
$\lambda_1$ and $\lambda_2$ can be switched.)

\begin{lem}  \label{norm-lem}
For $a,b,c \in \BB C$, we have:
\[  
\left\|\bigl(\begin{smallmatrix} a & 0 \\ c & 0 \end{smallmatrix}\bigr)\right\|
=\left\|\bigl(\begin{smallmatrix} 0 & a \\ 0 & c \end{smallmatrix}\bigr)\right\|
= \sqrt{|a|^2+|c|^2},
\]
\[  
\left\|\bigl(\begin{smallmatrix} a & b \\ 0 & 0 \end{smallmatrix}\bigr)\right\|
=\left\|\bigl(\begin{smallmatrix} 0 & 0 \\ a & b \end{smallmatrix}\bigr)\right\|
= \sqrt{|a|^2+|b|^2}.
\]
\end{lem}

\begin{proof}
For any $z= \bigl(\begin{smallmatrix} z_1 \\ z_2 \end{smallmatrix}\bigr)
\in \BB C^2$ with $\|z\|=1$, by the Cauchy-Schwarz inequality,
\[
\left\|\bigl(\begin{smallmatrix} a & b \\ 0 & 0 \end{smallmatrix}\bigr)
\bigl(\begin{smallmatrix} z_1 \\ z_2 \end{smallmatrix}\bigr)\right\|
= |az_1+bz_2| \le
\left\| \bigl(\begin{smallmatrix} a \\ b \end{smallmatrix}\bigr) \right\|
\cdot \left\| \bigl(\begin{smallmatrix} z_1 \\ z_2 \end{smallmatrix}\bigr)
\right\| = \sqrt{|a|^2+|b|^2}.
\]
On the other hand, if 
$\bigl(\begin{smallmatrix} z_1 \\ z_2 \end{smallmatrix}\bigr)
= (|a|^2+|b|^2)^{-\frac12} \bigl(\begin{smallmatrix} \bar a \\
  \bar b \end{smallmatrix}\bigr)$,
\[
\left\|\bigl(\begin{smallmatrix} a & b \\ 0 & 0 \end{smallmatrix}\bigr)
\bigl(\begin{smallmatrix} z_1 \\ z_2 \end{smallmatrix}\bigr)\right\|
= (|a|^2+|b|^2)^{-\frac12} |a \bar a + b \bar b|
= \sqrt{|a|^2+|b|^2}.
\]
This proves
$\left\|\bigl(\begin{smallmatrix} a & b \\ 0 & 0 \end{smallmatrix}\bigr)\right\|
= \sqrt{|a|^2+|b|^2}$.
The proofs of the other cases are similar.
\end{proof}
  
We define open domains of $\HC$
\begin{equation}  \label{D}
\BB D^+ = \{ Z \in \HC;\: \|Z\|<1 \}, \qquad
\BB D^- = \{ Z \in \HC;\: \|Z\|>1 \}.
\end{equation}
Alternatively, these domains can be expressed as
\[
\BB D^+ = \{ Z \in \HC;\: ZZ^*<1 \}, \qquad
\BB D^- = \{ Z \in \HC;\: ZZ^*>1 \},
\]
where the inequalities $ZZ^*<1$ and $ZZ^*>1$ mean that the self-adjoint matrix
$ZZ^*-1$ is negative and positive definite respectively.

We recall the definition of Shilov boundary in the context of several complex
variables.

\begin{df}
Let $\Omega$ be a bounded domain in $\BB C^n$, and let
${\cal A}(\Omega)$ be the algebra of functions continuous on the closure
$\overline{\Omega}$ and holomorphic inside $\Omega$.

A closed subset $S \subset \ \partial\Omega$ is a {\em boundary} for
${\cal A}(\Omega)$ if every function $f \in {\cal A}(\Omega)$ satisfies
\[
\max_{z \in \Omega} |f(z)| = \max_{z \in S} |f(z)|.
\]
Then the {\em Shilov boundary} of $\Omega$ is the unique, minimal closed
subset $S$ that satisfies the above property for all functions in
${\cal A}(\Omega)$.
\end{df}

\begin{prop}  \label{Shilov-prop}
The Shilov boundary of $\BB D^+$ is $U(2)$.
\end{prop}

\begin{proof}
Let $\overline{\BB D^+}$ denote the closure of $\BB D^+$, and let
$f: \overline{\BB D^+} \to \BB C$ be a continuous function that is
holomorphic on $\BB D^+$.
If $Z_0 \in \overline{\BB D^+}$ is such that $|f(Z_0)| \ge |f(Z)|$ for all
$Z \in \overline{\BB D^+}$, write $Z_0$ as a product $Z_0=U_0D_0V^*_0$
as in Proposition \ref{svd}.
Then the polydisk
\[
Pol = \left\{ U_0 \bigl(\begin{smallmatrix} z_1 & 0 \\
  0 & z_2 \end{smallmatrix}\bigr) V^*_0 ;\:
z_1, z_2 \in \BB C,\: |z_1|<1,\: |z_2|<1 \right\}
\]
is contained in $\BB D^+$ and $Z_0 \in Pol$.
Therefore, by the maximum modulus principle applied to $f$ restricted to
the closure of $Pol$,
\begin{multline*}
|f(Z_0)| \le \sup\{ |f(Z)|;\: Z \in \overline{Pol} \}
= \sup\{ |f(Z)|;\: Z \in \text{Shilov boundary of }Pol \}  \\
= \sup \left\{ |f(Z)|;\: Z = U_0
\bigl(\begin{smallmatrix} z_1 & 0 \\ 0 & z_2 \end{smallmatrix}\bigr) V^*_0 ;\:
z_1, z_2 \in \BB C,\: |z_1|=|z_2|=1 \right\}
\le  \sup\{ |f(Z)|;\: Z \in U(2) \}.
\end{multline*}
This proves that $U(2)$ contains the Shilov boundary of $\BB D^+$.

The other inclusion is proved by considering holomorphic functions of the form
$\tr(ZU^*+1)$, $U \in U(2)$.
\end{proof}

Equivalently, $U(2)$ consists of the boundary points for which both singular
values are equal to $1$. This description concerns singular values; the
ordinary eigenvalues of a unitary matrix may be distinct and play no role here.

For each $R>0$, we can rescale the open domains $\BB D^+$, $\BB D^-$
and define
\begin{align}
\BB D^+_R &= \{ Z \in \HC ;\: \|Z\|<R \} = \{ Z \in \HC ;\: ZZ^*<R^2 \},
\label{D^+_R}  \\
\BB D^-_R &= \{ Z \in \HC ;\: \|Z\|>R \} = \{ Z \in \HC ;\: ZZ^*>R^2 \}.
\label{D^-_R}
\end{align}
These domains $\BB D^+_R$ have
\[
U(2)_R = \{ Z \in \HC ;\: R^{-1}Z \in U(2) \}
\]
as their Shilov boundary.

\subsection{Conformal Transformations of $\HC$ Preserving $\BB D^+$}

We often regard the group $U(2,2)$ as a subgroup of $GL(2,\HC)$,
as described in Subsection 3.5 of \cite{FL1}. That is
\begin{multline}  \label{U(2,2)}
U(2,2) = \Bigl\{ \bigl(\begin{smallmatrix} a & b \\
  c & d \end{smallmatrix}\bigr) \in GL(2,\HC) ;\:
\bigl(\begin{smallmatrix} a & b \\ c & d \end{smallmatrix}\bigr)
\bigl(\begin{smallmatrix} 1 & 0 \\ 0 & -1 \end{smallmatrix}\bigr)
\bigl(\begin{smallmatrix} a^* & c^* \\ b^* & d^* \end{smallmatrix}\bigr)
= \bigl(\begin{smallmatrix} 1 & 0 \\ 0 & -1 \end{smallmatrix}\bigr) \Bigr\}  \\
= \Biggl\{ \begin{pmatrix} a & b \\ c & d \end{pmatrix}
\in GL(2,\HC) ;\: a,b,c,d \in \HC,\:
\begin{matrix} a^*a = 1+c^*c \\ d^*d = 1+b^*b \\ a^*b=c^*d \end{matrix}
\Biggr\}.
\end{multline}

\begin{prop}  \label{mobius-prop}
The group $U(2,2)$ acts on $\HC$ (with singularities) by fractional
linear transformations \eqref{conformal-action} preserving
$U(2) \subset \HC$ and open domains $\BB D^+$, $\BB D^-$.
\end{prop}

\begin{proof}
If
$\bigl(\begin{smallmatrix} a & b \\ c & d \end{smallmatrix}\bigr) \in U(2,2)$,
then
$\bigl(\begin{smallmatrix} a & b \\ c & d \end{smallmatrix}\bigr)^{-1}
= \bigl(\begin{smallmatrix} a^* & -c^* \\ -b^* & d^* \end{smallmatrix}\bigr)$.
The fractional linear transformation \eqref{conformal-action} associated to
$\bigl(\begin{smallmatrix} a & b \\ c & d \end{smallmatrix}\bigr)$
sends $Z \in \HC$ into $W = (a^*+Zb^*)^{-1}(c^*+Zd^*)$.
We have:
\begin{align*}
WW^*<1 \: &\Longleftrightarrow \:
(a^*+Zb^*)^{-1}(c^*+Zd^*)(c+dZ^*)(a+bZ^*)^{-1} <1  \\
&\Longleftrightarrow \:
(a^*+Zb^*)^{-1}\bigl( (c^*+Zd^*)(c+dZ^*) - (a^*+Zb^*)(a+bZ^*) \bigr)(a+bZ^*)^{-1}
<0  \\
&\Longleftrightarrow \:
(c^*+Zd^*)(c+dZ^*) - (a^*+Zb^*)(a+bZ^*) <0  \\
&\Longleftrightarrow \:
Z(d^*d-b^*b)Z^* - (a^*a-c^*c) + Z(d^*c-b^*a) + (c^*d-a^*b)Z^*  <0  \\
&\Longleftrightarrow \: ZZ^*<1.
\end{align*}
This proves that $U(2,2)$ preserves $\BB D^+$.
The proof for $\BB D^-$ is the same.
Then $U(2,2)$ preserves the boundary of $\BB D^+$ as well as its
Shilov boundary $U(2)$.
\end{proof}

Similarly, for each $R>0$ we can define a conjugate of $U(2,2)$
\[
U(2,2)_R = \bigl(\begin{smallmatrix} R & 0 \\ 0 & 1 \end{smallmatrix}\bigr)
U(2,2) \bigl(\begin{smallmatrix} R^{-1} & 0 \\ 0 & 1 \end{smallmatrix}\bigr)
\quad \subset GL(2,\HC).
\]
Each group $U(2,2)_R$ is a real form of $GL(2,\HC)$, preserves $U(2)_R$
and open domains $\BB D^+_R$ and $\BB D^-_R$.

Fix an element $Z_0 \in \BB D^+$ and express it as a product $Z_0=U_0D_0V^*_0$
as in Proposition \ref{svd} with $U_0, V_0 \in U(2)$ and
$D_0 = \bigl(\begin{smallmatrix} \lambda_1 & 0 \\
  0 & \lambda_2 \end{smallmatrix}\bigr)$,
$0 \le \lambda_1, \lambda_2 <1$. Write $(1-D_0^2)^{-\frac12}$ for the matrix
\[
\begin{pmatrix} \sqrt{1-\lambda_1^2} & 0 \\
  0 & \sqrt{1-\lambda_2^2} \end{pmatrix}^{-1} \in \HC.
\]

\begin{lem}  \label{Mobius-D-lem}
The matrix
\[
\begin{pmatrix} U_0 & U_0D_0 \\
  V_0D_0 & V_0 \end{pmatrix} (1-D_0^2)^{-\frac12}
= \begin{pmatrix} U_0 (1-D_0^2)^{-\frac12} & U_0 D_0(1-D_0^2)^{-\frac12} \\
  V_0 D_0(1-D_0^2)^{-\frac12} & V_0 (1-D_0^2)^{-\frac12} \end{pmatrix}
\]
belongs to $U(2,2)$, and the corresponding fractional linear transformation
\[
W \:\mapsto\:
U_0 (1-D_0^2)^{-\frac12} (W + D_0) (D_0W + 1)^{-1} (1-D_0^2)^{\frac12} V^*_0
\]
maps $0$ into $Z_0$.

Its inverse
\[
(1-D_0^2)^{-\frac12} \begin{pmatrix} U^*_0 & - D_0V^*_0 \\
  - D_0U^*_0 & V^*_0 \end{pmatrix}
= \begin{pmatrix} (1-D_0^2)^{-\frac12} U^*_0 & - D_0(1-D_0^2)^{-\frac12} V^*_0 \\
  - D_0(1-D_0^2)^{-\frac12} U^*_0 & (1-D_0^2)^{-\frac12} V^*_0 \end{pmatrix}
\]
also belongs to $U(2,2)$, and the corresponding fractional linear
transformation
\[
W \:\mapsto\:
(1-D_0^2)^{-\frac12} U^*_0 (W-Z_0) (1-Z_0^*W)^{-1} V_0 (1-D_0^2)^{\frac12}
\]
maps $Z_0$ into $0$.
\end{lem}

\begin{proof}
One checks that these matrices satisfy the conditions in \eqref{U(2,2)}.
The fact that the fractional linear transformations map $Z_0$ into $0$
and $0$ into $Z_0$ respectively is clear.
\end{proof}

\begin{cor}  \label{hom-space-cor}
  The group $U(2,2)$ acts on $\BB D^+$ transitively.
  Hence $\BB D^+$ can be naturally identified with the
  homogeneous space $U(2,2)/(U(2) \times U(2))$.
\end{cor}

\begin{proof}
Lemma \ref{Mobius-D-lem} implies that the action of $U(2,2)$ is transitive.
The stabilizer group of $0 \in \BB D^+$ in $U(2,2)$ is
\begin{equation}  \label{0-stabilizer}
\Bigl\{ \bigl(\begin{smallmatrix} a & 0 \\ 0 & d \end{smallmatrix}\bigr)
\in GL(2,\HC) ;\: a,d \in \HC,\:
\begin{smallmatrix} a^*a = 1 \\ d^*d = 1 \end{smallmatrix} \Bigr\}
\simeq U(2) \times U(2).
\end{equation}
\end{proof}

\subsection{Two Holomorphic Extension Lemmas}

In this subsection we prove two extension lemmas that can be regarded as
four-dimensional analogues of the following elementary fact about
holomorphic functions of one variable:
If a power series $\sum_{n=0}^{\infty} a_n z^n$ with coefficients
$a_n \in \BB C$ converges for all real $z \in (-1,1) \subset \BB R$,
then the series must converge on the entire unit disk
${\cal D} = \{ z \in \BB C ;\: |z| <1 \}$.

\begin{lem}  \label{non-zero-denom-lem}
Let $Z \in \BB D^+_R \cup \BB D^-_R$ and $W \in U(2)_R$.
Then $N(Z-W) \ne 0$.
\end{lem}

\begin{proof}
Write $Z = RV$ and $W=RU$ with $V \in \BB D^+ \cup \BB D^-$ and
$U \in U(2)$. Then
\[
N(Z-W)=R^2 \cdot N(V-U) = R^2 \cdot N(U) \cdot N(VU^{-1}-1) \ne 0,
\]
since $VU^{-1} \in \BB D^+ \cup \BB D^-$ and cannot have eigenvalue $1$.
\end{proof}

We have the following holomorphic extension lemmas.

\begin{lem}  \label{harm-ext-lem}
Let $\phi$ be a $\BB C$-valued harmonic function defined on the ball
${\cal B}_R = \{ X \in \BB H ;\: \|X\|_{\BB H}<R \}$.
Then $\phi$ has a holomorphic extension to a $\BB C$-valued
holomorphic harmonic function $\tilde \phi$ defined on $\BB D^+_R$.
\end{lem}

\begin{proof}
Suppose first that $\phi$ extends as a harmonic function to a neighborhood of
$\overline{{\cal B}_R}$.
The boundary of the ball $\partial{\cal B}_R$ lies in $U(2)_R$.
From the Poisson formula (Theorem \ref{Poisson-thm}) and
Lemma \ref{non-zero-denom-lem} we see that
\[
\tilde \phi(Z) = \frac1{2\pi^2} \int_{X \in \partial{\cal B}_R}
\frac{R^2-N(Z)}{N(X-Z)^2} \cdot \phi(X) \,\frac{dS}R,
\qquad Z \in \BB D^+,
\]
is a holomorphic extension of $\phi$ to a $\BB C$-valued holomorphic
harmonic function defined on $\BB D^+_R$.

In general, for each $r$, $0<r<R$, by above argument we can construct a
holomorphic extension $\tilde\phi_r$ of $\phi$ to a harmonic function
defined on $\BB D^+_r$. Moreover, this extension is unique, and
\[
\tilde\phi_{r_1}(Z) = \tilde\phi_{r_2}(Z) \qquad
\text{whenever $0<r_1<r_2$ and $Z \in \BB D^+_{r_1}$}.
\]
Thus, $\phi$ has a holomorphic extension to a holomorphic harmonic function
$\tilde \phi$ defined on the union $\bigcup_{0<r<R} \BB D^+_r = \BB D^+_R$.
\end{proof}

\begin{lem}  \label{reg-ext-lem}
Let $F$ be a left regular function defined on the ball
${\cal B}_R = \{ X \in \BB H ;\: \|X\|_{\BB H}<R \}$ with values in $\HC$ or
$\BB S$.
Then $F$ has a holomorphic extension to a holomorphic left regular function
$\tilde F$ defined on $\BB D^+_R$ with values in $\HC$ or $\BB S$ respectively.

Similarly, let $G$ be a right regular function defined on the ball ${\cal B}_R$
with values in $\HC$ or $\BB S'$.
Then $G$ has a holomorphic extension to a holomorphic right regular function
$\tilde G$ defined on $\BB D^+_R$ with values in $\HC$ or $\BB S'$ respectively.
\end{lem}

\begin{proof}
Suppose first that $F$ extends as a left regular function to a neighborhood of
$\overline{{\cal B}_R}$.
From the Cauchy-Fueter formula (Theorem \ref{Fueter}, observation
\eqref{Dx-restriction} and Lemma \ref{non-zero-denom-lem} we see that
\[
\tilde F(Z) = \frac1{2\pi^2} \int_{X \in \partial{\cal B}_R} \frac{(X-Z)^{-1}}{N(X-Z)}
\cdot \frac{X}{R} \cdot F(X) \,dS, \qquad Z \in \BB D^+.
\]
is a holomorphic extension of $\phi$ to a holomorphic left regular function
defined on $\BB D^+_R$ with values in $\HC$ or $\BB S$ respectively.

In general, as in the harmonic case, we extend $F$ to a holomorphic
left regular function $\tilde F_r$ defined on each $\BB D^+_r$, $0<r<R$,
then observe that such an extension is unique.
Hence, combining these extensions, we obtain a holomorphic extension of $F$
to a holomorphic left regular function $\tilde F$ defined on the union
$\bigcup_{0<r<R} \BB D^+_r = \BB D^+_R$.

The case of a right regular function is similar.
\end{proof}

For a fixed $Z_0 \in \HC$, the Cauchy-Fueter kernel
\[
\frac{(Z-Z_0)^{-1}}{N(Z-Z_0)}
\]
is both a holomorphic left regular and a holomorphic right regular function
wherever the denominator does not vanish.
Considering various values $Z_0 \in U(2)$, it is easy to show that $\BB D^+$
is maximal as a domain of left or right regular functions.
Thus, $\BB D^+$ can be thought of as a ``domain of regularity'' for
quaternionic regular functions.

\section{Schwarz Lemmas on $\BB D^+$}

In this section we prove a simple generalization of the classical
Schwarz Lemma to holomorphic functions on $\BB D^+$ with values in
a normed finite-dimensional complex vector space ${\cal V}$.
Examples of such normed vector spaces ${\cal V}$ include:
\begin{itemize}
\item
  Complex numbers $\BB C$ with the standard norm $|z|=\sqrt{x^2+y^2}$, $z=x+iy$.
\item
  Spinors $\BB S \simeq \BB C^2$ with norm
\[
  \left\| \bigl(\begin{smallmatrix} z_1 \\ z_2 \end{smallmatrix}\bigr) \right\|
  = \sqrt{|z_1|^2+|z_2|^2}.
\]
\item
  Spinors $\BB S' \simeq \BB C^2$ with norm
\[
  \| (z_1,\: z_2) \| = \sqrt{|z_1|^2+|z_2|^2}.
\]
\item
  Complexified quaternions $\HC$ regarded as complex $2 \times 2$ matrices
  and equipped with spectral norm \eqref{spectral-norm}.
\end{itemize}

\subsection{A Generalization of the Maximum Modulus Principle}

We start with a multivariable version of the maximum modulus principle.

\begin{thm}[Maximum Modulus Principle]  \label{MMP-thm}
Let $\Omega \subset \BB C^n$ be a non-empty connected open set, and
let ${\cal V}$ be a normed finite-dimensional complex vector space.
Suppose that $f: \Omega \to {\cal V}$ is holomorphic and that
$\|f(z)\| \le \|f(z_0)\|$ in a neighborhood of $z_0 \in \Omega$.
If $L:{\cal V}\to\BB C$ is a norm-one complex linear functional with
$L(f(z_0))=\|f(z_0)\|$, then $L\circ f$ is constant on $\Omega$.
In particular, if the norm on ${\cal V}$ is strictly convex, then $f$ is
constant.
\end{thm}

Indeed, $|L(f(z))|\le\|f(z)\|\le\|f(z_0)\|=L(f(z_0))$ near $z_0$, so the
scalar maximum modulus principle makes $L\circ f$ constant. If the norm is
strictly convex, the exposed face
$\{v;\ \|v\|\le\|f(z_0)\|,\ L(v)=\|f(z_0)\|\}$ consists only of $f(z_0)$,
so $f$ is locally, and hence globally, constant. For a general norm the last
conclusion need not hold; the scalarized statement above is the form used
below.

The corresponding boundary statement is:

\begin{thm}
Let $\Omega \subset \BB C^n$ be a bounded domain with closure
$\overline{\Omega}$, and let ${\cal V}$ be a normed finite-dimensional
complex vector space.
Suppose that $f: \overline{\Omega} \to {\cal V}$ is a continuous function that
is holomorphic on $\Omega$, then $\|f(z)\|$ attains a maximum at some point
of the Shilov boundary of $\Omega$. 
\end{thm}

\begin{proof}
Choose $z_1 \in \overline{\Omega}$ at which $\|f\|$ attains its
maximum. Applying the Hahn-Banach theorem to the complex vector space
${\cal V}$, and multiplying the functional by a unimodular scalar if
necessary, we can find a complex linear functional
$L: {\cal V} \to \BB C$ such that
\[
L(f(z_1))=\|f(z_1)\| \qquad \text{and} \qquad |L(v)| \le \|v\|
\text{ for all $v \in {\cal V}$}.
\]
Then the $\BB C$-valued function $L \circ f$ is holomorphic on $\Omega$ and
continuous on $\overline{\Omega}$.
Hence $|L(f(z))|$ attains a maximum at some point $z_0$
of the Shilov boundary of $\Omega$.
We have:
\[
\|f(z_0)\| \ge |L(f(z_0))| \ge |L(f(z_1))| = \|f(z_1)\| \ge \|f(z)\|
\qquad \text{for all $z \in \overline{\Omega}$}.
\]
This proves that $\|f(z)\|$ attains a maximum at some point
of the Shilov boundary of $\Omega$. 
\end{proof}

Combining this with Proposition \ref{Shilov-prop}, we obtain:

\begin{cor}  \label{MMP-cor}
Let ${\cal V}$ be a normed finite-dimensional complex vector space.
Suppose that $f: \overline{\BB D^+} \to {\cal V}$ is a continuous function
that is holomorphic on $\BB D^+$, then $\|f(z)\|$ attains a maximum at
some point of $U(2)$. 
\end{cor}

\subsection{A Generalization of the Schwarz Lemma}  \label{Schwarz-subsect}

If $f$ is holomorphic in a neighborhood of a point $a \in \BB C^n$, we say
that {\em $f$ is of total order $k$ at $a$} if all the partial derivatives
at $a$ of orders less than $k$ vanish:
\[
\frac{\partial^{k_1+\ldots+k_n} f}{\partial z_1^{k_1} \dots \partial z_n^{k_n}}(a)=0
\quad \text{for all $k_1+\cdots+k_n<k$}.
\]
The partial derivative of order $0$ is, of course, $f$ itself.
Thus, $f$ being of total order $1$ at $a$ amounts to $f(a)=0$.

Recall that $\overline{\BB D^+}$ denotes the closure of $\BB D^+$.

\begin{thm}[Schwarz Lemma]
Let ${\cal V}$ be a normed finite-dimensional complex vector space, and
let $f$ be a ${\cal V}$-valued function holomorphic on a neighborhood of
$\overline{\BB D^+}$, and suppose that $f$ is of total order $k$ at $0$
and that $\|f(Z)\| \le M$ for all $Z \in \overline{\BB D^+}$. Then
\begin{equation}  \label{SLem-eqn}
\|f(Z)\| \le M \cdot \|Z\|^k
\qquad \text{ for all $Z \in \overline{\BB D^+}$}.
\end{equation}
\end{thm}

\begin{proof}
This proof is adapted from the proof of Schwarz Lemma for polydisks
given in \cite{GR} (Theorem 8).
Let $Z \in \BB D^+$, $Z \ne 0$.
Observe that $tZ/\|Z\| \in \BB D^+$ for all
$t \in {\cal D} = \{ z \in \BB C;\: |t|<1 \}$.
Expanding $f$ into homogeneous polynomials, we see that
\[
g(t) = t^{-k} \cdot f(tZ/\|Z\|)
\]
is a function holomorphic in ${\cal D}$ (no singularity at $t=0$)
and $\|g(t)\| \le M$ for $|t|=1$.
Thus, by the maximum modulus principle, $\|g(t)\| \le M$ for all $t$,
$|t| \le 1$.
Setting $t=\|Z\|$, it follows that
\[
\|Z\|^{-k} \cdot \| f(Z) \| = \bigl\| \|Z\|^{-k} \cdot f(Z) \bigr\|
= \| g(\|Z\|) \| \le M,
\]
and $\|f(Z)\| \le M \cdot \|Z\|^k$.
Since this inequality holds for all $Z \in \BB D^+$, \eqref{SLem-eqn} follows.
\end{proof}

\begin{cor}  \label{Schwarz-cor0}
Let ${\cal V}$ be a normed finite-dimensional complex vector space, and
let $f: \BB D^+ \to {\cal V}$ be a holomorphic function such that
$f(0)=0$ and $\|f(Z)\| \le M$ for all $Z \in \BB D^+$, then
\[
\|f(Z)\| \le M \cdot \|Z\| \qquad \text{for all $Z \in \BB D^+$},
\]
and the partial derivatives at the origin satisfy
\[
\biggl\| \frac{\partial f}{\partial z_{ij}}(0) \biggr\| \le M,
\qquad i,j \in \{1,2\}.
\]
\end{cor}

\begin{cor}  \label{Schwarz-cor1}
Let $f: \BB D^+ \to \BB C$ be a holomorphic function such that $f(0)=0$
and $|f(Z)| \le M$ for all $Z \in \BB D^+$, then
\[
|f(Z)| \le M \cdot \|Z\| \qquad \text{for all $Z \in \BB D^+$},
\]
and the partial derivatives at the origin satisfy
\[
\biggl| \frac{\partial f}{\partial z_{ij}}(0) \biggr| \le M,
\qquad i,j \in \{1,2\}.
\]
\end{cor}

\begin{cor}  \label{Schwarz-cor2}
Let $F$ be a holomorphic left regular function on $\BB D^+$ with values in
$\BB S$ or $\HC$.
Suppose that $F(0)=0$ and $\|F(Z)\| \le M$ for all $Z \in \BB D^+$, then
\[
\|F(Z)\| \le M \cdot \|Z\| \qquad \text{for all $Z \in \BB D^+$},
\]
and the partial derivatives at the origin satisfy
\[
\biggl\| \frac{\partial F}{\partial z_{ij}}(0) \biggr\| \le M,
\qquad i,j \in \{1,2\}.
\]

Similarly, Let $G$ be a holomorphic right regular function on $\BB D^+$ with
values in $\BB S'$ or $\HC$.
Suppose that $G(0)=0$ and $\|G(Z)\| \le M$ for all $Z \in \BB D^+$, then
\[
\|G(Z)\| \le M \cdot \|Z\| \qquad \text{for all $Z \in \BB D^+$},
\]
and the partial derivatives at the origin satisfy
\[
\biggl\| \frac{\partial G}{\partial z_{ij}}(0) \biggr\| \le M,
\qquad i,j \in \{1,2\}.
\]
\end{cor}

\subsection{Counterexamples}

The second part of the classical Schwarz Lemma asserts that if
$f: {\cal D} \to {\cal D}$ is a holomorphic function with $f(0)=0$ and
$|f(z_0)|=|z_0|$ for some $z_0 \in {\cal D} \setminus \{0\}$,
then $f$ is of the form $f(z) = \lambda z$, where $|\lambda|=1$.
However, an analogous statement for holomorphic functions on $\BB D^+$
fails, and in this subsection we give two counterexamples -- one involving
a harmonic function and the other involving a left regular function.

We consider a holomorphic harmonic function $\phi: \HC \to \BB C$
\[
\phi(Z) = \phi \bigl(\begin{smallmatrix} z_{11} & z_{12} \\
  z_{21} & z_{22} \end{smallmatrix}\bigr) = z_{11} + \epsilon (z_{12})^2,
\]
where $\epsilon$ is a positive real parameter.
Clearly, $\phi(0)=0$, $\phi$ is {\em not} linear, and if $Z_0 \in \BB D^+$
is of the form
$\bigl(\begin{smallmatrix} a & 0 \\ 0 & 0 \end{smallmatrix}\bigr)$,
then, by Lemma \ref{norm-lem},  $|\phi(Z_0)|=\|Z_0\|$.
On the other hand, for $0 < \epsilon < 1/2$, we can show that
$\sup\{ |\phi(Z)|;\: Z \in \overline{\BB D^+} \} =1$.
Indeed, $|\phi(1)|=1$ and, by Corollary \ref{MMP-cor},
\[
\sup\{ |\phi(Z)|;\: Z \in \overline{\BB D^+} \}
= \sup\{ |\phi(Z)|;\: Z \in U(2) \}.
\]
When $Z \in U(2)$, $Z= \bigl(\begin{smallmatrix} z_{11} & z_{12} \\
  z_{21} & z_{22} \end{smallmatrix}\bigr)$ with $|z_{11}|^2+|z_{12}|^2=1$.
Hence
\[
|\phi(Z)| \le |z_{11}|+\epsilon|z_{12}|^2
= \epsilon + |z_{11}| - \epsilon|z_{11}|^2 \le 1.
\]

Next, we consider a holomorphic left regular function $F: \HC \to \BB S$
\[
F(Z) = F \bigl(\begin{smallmatrix} z_{11} & z_{12} \\
  z_{21} & z_{22} \end{smallmatrix}\bigr)
= \begin{pmatrix} z_{11} + \epsilon (z_{12})^2 \\
  2\epsilon z_{11} z_{12} \end{pmatrix},
\]
where $\epsilon$ is a positive real parameter.
As in the previous example, $F(0)=0$, $F$ is {\em not} linear, and if
$Z_0 \in \BB D^+$ is of the form
$\bigl(\begin{smallmatrix} a & 0 \\ 0 & 0 \end{smallmatrix}\bigr)$,
then, by Lemma \ref{norm-lem},  $\|F(Z_0)\|=\|Z_0\|$.
On the other hand, for $0 < \epsilon < 1/4$, we can show that
$\sup\{ \|F(Z)\|;\: Z \in \overline{\BB D^+} \} =1$.
Indeed, $\|F(1)\|=1$ and, by Corollary \ref{MMP-cor},
\[
\sup\{ \|F(Z)\|;\: Z \in \overline{\BB D^+} \}
= \sup\{ \|F(Z)\|;\: Z \in U(2) \}.
\]
As before, when $Z \in U(2)$,
$Z= \bigl(\begin{smallmatrix} z_{11} & z_{12} \\
  z_{21} & z_{22} \end{smallmatrix}\bigr)$
with $|z_{11}|^2+|z_{12}|^2=1$. Hence
\[
\|F(Z)\|^2 =
\bigl| z_{11} + \epsilon (z_{12})^2 \bigr|^2 + 4\epsilon^2|z_{11}|^2|z_{12}|^2
\le \bigl(\epsilon + |z_{11}| - \epsilon|z_{11}|^2\bigr)^2
+ 4\epsilon^2|z_{11}|^2(1-|z_{11}|^2).
\]
Let $t=|z_{11}|$, $t \in [0,1]$. Then
$\|F(Z)\|^2 \le (\epsilon + t - \epsilon t^2)^2 + 4\epsilon^2 t^2(1-t^2)$.
Differentiating with respect to $t$, we see that
\begin{multline*}
\bigl( (\epsilon + t - \epsilon t^2)^2 + 4\epsilon^2 t^2(1-t^2) \bigr)'
= 2(1-2\epsilon t)(\epsilon + t - \epsilon t^2) + 8\epsilon^2t(1-2t^2)  \\
\ge 2(1-2\epsilon t)t - 8\epsilon^2t
\ge 2t(1-2\epsilon-4\epsilon^2) \ge 0.
\end{multline*}
Thus, $h(t)= (\epsilon + t - \epsilon t^2)^2 + 4\epsilon^2 t^2(1-t^2)$
is strictly increasing on $[0,1]$ with maximal value $h(1)=1$.
This shows that $\|F(Z)\|^2 \le 1$ for all $Z \in U(2)$
and hence all $Z \in \overline{\BB D^+}$.

A similar counterexample can be constructed using holomorphic right regular
functions.

\subsection{A Generalization of the Schwarz--Pick Theorem}  \label{Schwarz-Pick-subsect}

In this subsection we combine our generalizations of the Schwarz Lemma
with fractional linear transformations of $\BB D^+$ and ${\cal D}$ to obtain
two analogues of the Schwarz--Pick Theorem.
First, we have an analogue for $\BB C$-valued functions.

\begin{thm}
Let $f: \BB D^+ \to {\cal D}$ be holomorphic, and let $Z_1 \in \BB D^+$.
Express $Z_1=U_1D_1V^*_1$ as in Proposition \ref{svd}.
Then, for all $Z_2 \in \BB D^+$, $Z_2 \ne Z_1$, we have
\begin{align*}
\left| \frac{f(Z_1)-f(Z_2)}{1-\overline{f(Z_1)}f(Z_2)}\right|
&\le \left\| (1-D_1^2)^{-\frac12} U^*_1 (Z_1-Z_2)
(1-Z_1^*Z_2)^{-1} V_1 (1-D_1^2)^{\frac12} \right\|  \\
&\le \left\| (1-D_1^2)^{\frac12} \right\|
\cdot \left\| (1-D_1^2)^{-\frac12} \right\| \cdot
\left\| (Z_1 - Z_2) (1-Z_1^* Z_2)^{-1} \right\|.
\end{align*}
\end{thm}

\begin{proof}
Consider fractional linear transformations $\phi: {\cal D} \to {\cal D}$
defined by
\[
\phi(\zeta)= \frac{\zeta-f(Z_1)}{1-\overline{f(Z_1)}\zeta}
\]
and $\psi: \BB D^+ \to \BB D^+$ defined by
\[
\psi(Z)=
U_1 (1-D_1^2)^{-\frac12} (Z + D_1) (D_1Z + 1)^{-1} (1-D_1^2)^{\frac12} V^*_1
\]
(see Lemma \ref{Mobius-D-lem}).
Then the composition $g = \phi \circ f \circ \psi$ is a holomorphic function
$\BB D^+ \to {\cal D}$ such that $g(0)=0$.
By Corollary \ref{Schwarz-cor1}, $|g(Z)| \le \|Z\|$ for all $Z \in \BB D^+$.
Substituting $Z=\psi^{-1}(Z_2)$, we obtain:
\begin{multline*}
\left| \frac{f(Z_1)-f(Z_2)}{1-\overline{f(Z_1)}f(Z_2)}\right|
= \bigl| g\bigl( \psi^{-1}(Z_2)\bigr)\bigr|
\le \|\psi^{-1}(Z_2)\|  \\
= \left\| (1-D_1^2)^{-\frac12} U^*_1 (Z_2-Z_1)
(1-Z_1^*Z_2)^{-1} V_1 (1-D_1^2)^{\frac12} \right\|  \\
\le \left\| (1-D_1^2)^{-\frac12} \right\| \cdot
\left\| (Z_2 - Z_1)(1 - Z_1^* Z_2)^{-1} \right\| \cdot
\left\| (1-D_1^2)^{\frac12} \right\|,
\end{multline*}
since $\|AB\| \le \|A\| \cdot \|B\|$ and $\|U^*_1 A V_1\|=\|A\|$.
\end{proof}

And here is an analogue of the Schwarz--Pick Theorem for $\HC$-valued functions.

\begin{thm}  \label{Schwarz-Pick-thm2}
Let $F: \BB D^+ \to \BB D^+$ be holomorphic, and let $Z_1 \in \BB D^+$.
Express $Z_1=U_1D_1V^*_1$ and $F(Z_1) = \tilde U_1 \tilde D_1 \tilde V^*_1$
as in Proposition \ref{svd}.
Then, for all $Z_2 \in \BB D^+$, $Z_2 \ne Z_1$, we have
\begin{multline*}
\left\| (1-\tilde D_1^2)^{-\frac12} \tilde U^*_1 \bigl( F(Z_1) - F(Z_2) \bigr)
\bigl(1-F(Z_1)^*F(Z_2)\bigr)^{-1} \tilde V_1(1-\tilde D_1^2)^{\frac12} \right\| \\ 
\le \left\| (1-D_1^2)^{-\frac12} U^*_1 (Z_1-Z_2)
(1-Z_1^*Z_2)^{-1} V_1 (1-D_1^2)^{\frac12} \right\|  \\
\le \left\| (1-D_1^2)^{\frac12} \right\|
\cdot \left\| (1-D_1^2)^{-\frac12} \right\| \cdot
\left\| (Z_1 - Z_2) (1-Z_1^* Z_2)^{-1} \right\|.
\end{multline*}
\end{thm}

\begin{proof}
Consider fractional linear transformations $\psi$ and
$\tilde\psi: \BB D^+ \to \BB D^+$ defined by
\begin{align*}
\psi(Z) &=
U_1 (1-D_1^2)^{-\frac12} (Z + D_1) (D_1Z + 1)^{-1} (1-D_1^2)^{\frac12} V^*_1,  \\
\tilde\psi(W) &=
(1-\tilde D_1^2)^{-\frac12} \tilde U^*_1 \bigl( W  - F(Z_1) \bigr)
\bigl( 1-F(Z_1)^* W \bigr)^{-1} \tilde V_1 (1-\tilde D_1^2)^{\frac12}.
\end{align*}
(see Lemma \ref{Mobius-D-lem}).
Then the composition $G = \tilde\psi \circ F \circ \psi$ is a holomorphic
function $\BB D^+ \to \BB D^+$ such that $G(0)=0$.
By Corollary \ref{Schwarz-cor0}, $\|G(Z)\| \le \|Z\|$ for all $Z \in \BB D^+$.
Substituting $Z=\psi^{-1}(Z_2)$, we obtain:
\begin{multline*}
\left\| (1-\tilde D_1^2)^{-\frac12} \tilde U^*_1 \bigl( F(Z_1) - F(Z_2) \bigr)
\bigl(1-F(Z_1)^*F(Z_2)\bigr)^{-1} \tilde V_1(1-\tilde D_1^2)^{\frac12} \right\|
= \bigl\| G\bigl( \psi^{-1}(Z_2)\bigr)\bigr\|  \\
\le \|\psi^{-1}(Z_2)\|
= \left\| (1-D_1^2)^{-\frac12} U^*_1 (Z_2-Z_1)
(1-Z_1^*Z_2)^{-1} V_1 (1-D_1^2)^{\frac12} \right\|  \\
\le \left\| (1-D_1^2)^{-\frac12} \right\| \cdot
\left\| (Z_2 - Z_1)(Z_1^* Z_2 - 1)^{-1} \right\| \cdot
\left\| (1-D_1^2)^{\frac12} \right\|,
\end{multline*}
since $\|AB\| \le \|A\| \cdot \|B\|$ and $\|U^*_1 A V_1\|=\|A\|$.
\end{proof}

The inequality in Theorem \ref{Schwarz-Pick-thm2} simplifies when
the matrix $\tilde D_1$ is scalar.

\begin{cor}
Let $F: \BB D^+ \to \BB D^+$ be holomorphic, and let $Z_1 \in \BB D^+$.
Suppose that $F(Z_1) = \tilde r \tilde U$ for some $\tilde r \in \BB R$
and $\tilde U \in U(2)$.
Express $Z_1=U_1D_1V^*_1$ as in Proposition \ref{svd}.
Then, for all $Z_2 \in \BB D^+$, $Z_2 \ne Z_1$, we have
\begin{multline*}
\left\|\bigl(F(Z_1)-F(Z_2)\bigr) \bigl(1-F(Z_1)^*F(Z_2)\bigr)^{-1} \right\| \\ 
\le \left\| (1-D_1^2)^{-\frac12} U^*_1 (Z_1-Z_2)
(1-Z_1^*Z_2)^{-1} V_1 (1-D_1^2)^{\frac12} \right\|  \\
\le \left\| (1-D_1^2)^{\frac12} \right\|
\cdot \left\| (1-D_1^2)^{-\frac12} \right\| \cdot
\left\| (Z_1 - Z_2) (1-Z_1^* Z_2)^{-1} \right\|.
\end{multline*}
\end{cor}

\begin{proof}
The assumption that $F(Z_1) = \tilde r \tilde U$ for some
$\tilde r \in \BB R$ and $\tilde U \in U(2)$ is equivalent to the matrix
$\tilde D_1$ in Theorem \ref{Schwarz-Pick-thm2} being scalar.
\end{proof}

\begin{rem}
Each of the above versions of the Schwarz--Pick Theorem simplifies further when
the matrix $D_1$ is scalar.
This is equivalent to $Z_1 \in \BB D^+$ being of the form 
$Z_1 = rU$ for some $r \in \BB R$ with $|r|<1$ and $U \in U(2)$.
\end{rem}

\begin{rem}
Theorem \ref{Schwarz-Pick-thm2} can also be applied to functions with values in
$\BB S$ and $\BB S'$ (such as left and right regular functions)
via norm-preserving inclusions
\[
\BB S \hookrightarrow \HC, \quad
\left(\begin{smallmatrix} z_1 \\ z_2 \end{smallmatrix}\right)
\mapsto \left(\begin{smallmatrix} z_1 & 0 \\ z_2 & 0 \end{smallmatrix}\right),
\qquad \text{and} \qquad 
\BB S' \hookrightarrow \HC, \quad (z_1,\: z_2)
\mapsto \left(\begin{smallmatrix} z_1 & z_2 \\ 0 & 0 \end{smallmatrix}\right),
\]
(Lemma \ref{norm-lem}).
\end{rem}

\section{The Poincar\'e Metric on $\BB D^+$}  \label{Poincare-sect}

Let $Z \in \BB D^+$, and express it as a product $Z=UDV^*$ as in
Proposition \ref{svd}.
The tangent space $T_Z\BB D^+$ is naturally isomorphic to $\HC$.
For an element ${\bf v} \in \HC$, we denote by $\|{\bf v}\|_{HS}$
the Hilbert--Schmidt norm of ${\bf v}$:
\[
\|{\bf v}\|_{HS} = \sqrt{\tr({\bf v}^*{\bf v})} =
\sqrt{|v_{11}|^2+|v_{12}|^2+|v_{21}|^2+|v_{22}|^2}, \qquad
{\bf v} = \bigl(\begin{smallmatrix} v_{11} & v_{12} \\
  v_{21} & v_{22} \end{smallmatrix} \bigr) \in \HC.
\]
Then we define the {\em Poincar\'e metric} on $\BB D^+$ by declaring that,
for each tangent vector ${\bf v} \in T_Z\BB D^+ \simeq \HC$,
\begin{equation}  \label{|v|_Z-def}
|{\bf v}|_Z = \bigl\| (1-D^2)^{-\frac12} U^*{\bf v}V
(1-D^2)^{-\frac12} \bigr\|_{HS}.
\end{equation}
This does not depend on the chosen singular value decomposition: equivalently,
by Lemma \ref{intrinsic-metric-lem} below, it is
$\|(1-ZZ^*)^{-1/2}{\bf v}(1-Z^*Z)^{-1/2}\|_{HS}$.

\begin{prop}  \label{Poincere-invar-prop}
The Poincar\'e metric on $\BB D^+$ is invariant under the action
of $U(2,2)$ on $\BB D^+$ by fractional linear transformations
\eqref{conformal-action}.
\end{prop}

\begin{proof}
The stabilizer of $0 \in \BB D^+$ in $U(2,2)$ can be identified with
$U(2) \times U(2)$ as in \eqref{0-stabilizer}.
It is easy to see that $U(2) \times U(2)$ preserves the Poincar\'e metric.
And it is sufficient to prove that, for any $g \in U(2,2)$, the tangent map
$g_*: T_0\BB D^+ \to T_{g(0)}\BB D^+$ preserves the Poincar\'e metric.
Express $g(0)$ as a product $g(0)=UDV^*$ as in Proposition \ref{svd}.
Then $g$ can be written as a product $g= \tilde g \cdot u$,
where $u \in U(2) \times U(2)$ and
$\tilde g=\bigl(\begin{smallmatrix} U & UD \\
  VD & V \end{smallmatrix}\bigr) (1-D^2)^{-\frac12}$ is a matrix from
Lemma \ref{Mobius-D-lem}.
Since $g_*=\tilde g_* \circ u_*$, it is sufficient to show that
$\tilde g_*: T_0\BB D^+ \to T_{g(0)}\BB D^+$ preserves the Poincar\'e metric.
By direct computation, $\tilde g_*$ maps
\[
T_0\BB D^+ \ni \quad {\bf v} \:\mapsto\:
U(1-D^2)^{\frac12} {\bf v} (1-D^2)^{\frac12} V^* \quad \in T_{g(0)}\BB D^+.
\]
The result follows from the invariance of the Hilbert--Schmidt norm under
multiplication by unitary matrices.
\end{proof}

\begin{prop}  \label{Poincare-dist-prop}
Let $Z \in \BB D^+$, and express it as a product $Z=UDV^*$ as in
Proposition \ref{svd} with
$D = \bigl(\begin{smallmatrix} \lambda_1 & 0 \\
  0 & \lambda_2 \end{smallmatrix}\bigr)$,
$0 \le \lambda_1, \lambda_2 < 1$, as in \eqref{Diag}.
The Poincar\'e distance of $0$ to $Z$ is
\begin{equation}  \label{dist-formula}
d(0,Z)= \frac12 \sqrt{ \Bigl(\log\frac{1+\lambda_1}{1-\lambda_1}\Bigr)^2 +
\Bigl(\log\frac{1+\lambda_2}{1-\lambda_2}\Bigr)^2}.
\end{equation}
\end{prop}

\begin{proof}
By the $U(2) \times U(2)$ invariance, we may assume that
$Z = D = \bigl(\begin{smallmatrix} \lambda_1 & 0 \\
  0 & \lambda_2 \end{smallmatrix}\bigr)$.
First, we show that the distance from $0$ to $D$ can be computed using curves
of the form
$\gamma(t)= \bigl(\begin{smallmatrix} \lambda_1(t) & 0 \\
  0 & \lambda_2(t) \end{smallmatrix}\bigr)$, $t \in [0,1]$, with $\gamma(0)=0$
and $\gamma(1)=D$. After subdividing at singular-value crossings, a
piecewise ${\cal C}^1$ curve from $0$ to $D$ inside $\BB D^+$ can be
written as $\gamma(t)= U(t) \cdot D(t) \cdot V^*(t)$ with
$U(t), V(t) \in U(2)$ and $D(t)$ diagonal as in Proposition \ref{svd}.
Observe that $U'(t) \cdot D(t) \cdot V^*(t)$ and $U(t) \cdot D(t) \cdot V'^*(t)$
are orthogonal to $U(t) \cdot D'(t) \cdot V^*(t)$ with respect to the
Poincar\'e metric. Hence
\begin{multline*}
\operatorname{Length}(\gamma(t)) = \int_0^1 |\gamma'(t)|_{\gamma(t)} \,dt  \\
= \int_0^1 \sqrt{|U(t) \cdot D'(t) \cdot V^*(t)|_{\gamma(t)}^2 +
|U'(t) \cdot D(t) \cdot V^*(t)+U(t) \cdot D(t) \cdot V'^*(t)|_{\gamma(t)}^2}\,dt\\
\ge \int_0^1 |U(t) \cdot D'(t) \cdot V^*(t)|_{\gamma(t)} \,dt
= \int_0^1 |D'(t)|_{D(t)} \,dt
= \operatorname{Length}(D(t)),
\end{multline*}
and $D(t)$ is a curve from $0$ to $D$ of the form
$D(t)= \bigl(\begin{smallmatrix} \lambda_1(t) & 0 \\
  0 & \lambda_2(t) \end{smallmatrix}\bigr)$.

Thus, the problem reduces to finding the distance between $0$ and
$D=(\lambda_1,\lambda_2)$ inside the open square $S = (-1,1) \times (-1,1)$
with Riemannian metric
\[
ds^2 = \frac{dx \otimes dx}{(1-x^2)^2} + \frac{dy \otimes dy}{(1-y^2)^2}.
\]
It is easy to see that $(S,ds^2)$ is isometric to $\BB R^2$ with the standard
Euclidean metric via the map
\[
(x,y) \mapsto \biggl( \frac12 \log\frac{1+x}{1-x},
\frac12 \log\frac{1+y}{1-y} \biggr).
\]
This map sends $0 \in S$ into the origin and
$D = \bigl(\begin{smallmatrix} \lambda_1 & 0 \\
  0 & \lambda_2 \end{smallmatrix}\bigr)\in S$ into
$\bigl( \frac12 \log\frac{1+\lambda_1}{1-\lambda_1},
\frac12 \log\frac{1+\lambda_2}{1-\lambda_2} \bigr)$,
and \eqref{dist-formula} follows.
\end{proof}
  
\begin{cor}
Suppose that $Z \in \BB D^+$ is of the form $Z = rU$ for some $r \in \BB R$
with $0 \le r <1$ and $U \in U(2)$.
Then the Poincar\'e distance of $0$ to $Z$ is
\[
d(0,Z) = \frac1{\sqrt{2}} \log \frac{1+r}{1-r}
= \frac1{\sqrt{2}} \log \frac{1+\|Z\|}{1-\|Z\|}.
\]
\end{cor}

\begin{proof}
The assumption that $Z = rU$ for some $r \ge 0$ and $U \in U(2)$
is equivalent to the matrix $D$ in Proposition \ref{Poincare-dist-prop}
being scalar:
$D=\bigl(\begin{smallmatrix} r & 0 \\ 0 & r \end{smallmatrix}\bigr)$.
\end{proof}

We address the uniqueness of the Poincar\'e metric.

\begin{prop}  \label{P-uniqueness-prop}
If a Riemannian metric on $\BB D^+$ is invariant under the action of
$U(2,2)$ on $\BB D^+$ by fractional linear transformations
\eqref{conformal-action}, then this metric is a constant multiple of the
Poincar\'e metric.
\end{prop}

\begin{proof}
Since $U(2,2)$ acts on $\BB D^+$ transitively, it is sufficient to show that
the two metrics are proportional on the tangent space $T_0 \BB D^+$.
Recall that the stabilizer of $0 \in \BB D^+$ can be identified with
$U(2) \times U(2)$ as in \eqref{0-stabilizer}.

\begin{lem}  \label{real-Schur-lem}
As a representation of $U(2) \times U(2)$ over $\BB R$, the tangent space
$T_0 \BB D^+$ is irreducible, and the algebra of endomorphisms intertwining
the $U(2) \times U(2)$ action
\[
\End_{U(2) \times U(2)} (T_0 \BB D^+)
\]
is isomorphic to $\BB C$.
\end{lem}

\begin{proof}
As a representation of $U(2) \times U(2)$ over $\BB C$, $T_0 \BB D^+$ is
irreducible, since it is the tensor product of two irreducible 2-dimensional
representations of $U(2)$. Restricting the scalars from $\BB C$ to $\BB R$
preserves the irreducibility because the multiplication by $i$ can be
realized by action of $U(2) \times U(2)$.

Comparing the actions of the center of $U(2) \times U(2)$, we see that,
as a representation over $\BB C$, $T_0 \BB D^+$ is not isomorphic to the
conjugate representation. By Theorems 6.3, 6.7 of Chapter II in \cite{BD},
$\End_{U(2) \times U(2)} (T_0 \BB D^+)$ is isomorphic to $\BB C$.
\end{proof}

We continue the proof of Proposition \ref{P-uniqueness-prop}.
The Poincar\'e metric produces a symmetric non-degenerate
$U(2) \times U(2)$-invariant bilinear form on $T_0 \BB D^+$,
denote it by $B({\bf v_1}, {\bf v_2})$.
Any other non-degenerate $U(2) \times U(2)$-invariant bilinear form on
$T_0 \BB D^+$ (that is not necessarily symmetric) can be obtained from $B$
as $B(A{\bf v_1}, {\bf v_2})$, where $A$ is an invertible element of
$\End_{U(2) \times U(2)} (T_0 \BB D^+)$. By Lemma \ref{real-Schur-lem},
$\End_{U(2) \times U(2)} (T_0 \BB D^+) \simeq \BB C$.
The multiplication by a $\lambda \in \BB C$ results in a bilinear form
$B(\lambda {\bf v_1}, {\bf v_2})$ which is symmetric if and only if
$\lambda$ is a real scalar.
Hence any other symmetric non-degenerate $U(2) \times U(2)$-invariant metric
on $T_0 \BB D^+$ is proportional to the Poincar\'e metric.
\end{proof}


\section{Half Plane Realizations of $\BB D^+$ and $\BB D^-$}  \label{half-plane-sect}

In this section we describe realizations of $\BB D^+$ and $\BB D^-$ that
are analogous to half plane realizations of the unit disk ${\cal D}$.
Then we describe the Poincar\'e metric on the analogues of half planes
$\BB T^+$ and $\BB T^-$.

We recall some results from Section 3 in \cite{FL1}.
Thus, the Minkowski space $\BB M$ is realized inside $\HC$ as
\begin{align*}
\BB M &= \bigl\{ -i{\bf e_0}y_0 + {\bf e_1}y_1 + {\bf e_2}y_2 + {\bf e_3}y_3 =
\bigl(\begin{smallmatrix} -iy_0-iy_3 & -iy_1-y_2 \\
  -iy_1+y_2 & -iy_0+iy_3 \end{smallmatrix}\bigr)
;\: y_0, y_1, y_2, y_3 \in \BB R \bigr\}  \\
&= \{Z \in \HC ; \: Z^* = - Z \}
= \bigl\{ Z= \bigl(\begin{smallmatrix} z_{11} & z_{12} \\
  z_{21} & z_{22} \end{smallmatrix}\bigr) \in \HC
; \: z_{11}, z_{22} \in i \BB R, \: z_{21} = -\overline{z_{12}} \bigr\}.
\end{align*}
Let $\C^+$ denote the open cone
\begin{equation}  \label{C^+}
\C^+ = \{ Y \in \BB M ;\: N(Y)<0, \: i\tr Y <0 \},
\end{equation}
and define two open tube domains in $\HC$
\[
\BB T^+ = \BB M + i\C^+,
\qquad \BB T^- = \BB M - i\C^+.
\]
Note that $i\C^+$ (respectively $-i\C^+$) consists of self-adjoint elements of
$\HC$ that are negative definite (respectively positive definite).
Recall a quaternionic analogue of the Cayley transform.

\begin{lem} [Lemmas 54 and 63 in \cite{FL1}]
Consider an element
$\gamma = \frac1{\sqrt{2}} \bigl(\begin{smallmatrix} -i & -i \\
  1 & -1 \end{smallmatrix}\bigr) \in GL(2,\HC)$
with
$\gamma^{-1} = \frac1{\sqrt{2}} \bigl(\begin{smallmatrix} i & 1 \\
  i & -1 \end{smallmatrix}\bigr)$.
The associated fractional linear map on $\HC$
\[
\gamma: \: Z \mapsto -i(Z+1)(Z-1)^{-1}
\]
maps $\BB M \to U(2)$, $\BB T^+ \to \BB D^+$ and $\BB T^- \to \BB D^-$.
This map has singularities on $\BB T^-$ and no singularities on
$\BB M \cup \BB T^+$.

Conversely, the fractional linear map on $\HC$
\[
\gamma^{-1}: \: Z \mapsto (Z-i)(Z+i)^{-1}
\]
maps $U(2) \to \BB M$, $\BB D^+ \to \BB T^+$ and $\BB D^- \to \BB T^-$.
This map has singularities on $U(2)$ and no singularities on
$\BB D^+ \cup \BB D^-$.
\end{lem} 

\begin{proof}
Let $Z \in \BB M$, so that $Z^*=-Z$, and let $W=-i(Z+1)(Z-1)^{-1}$.
We compute:
\[
W^*= i(Z^*-1)^{-1}(Z^*+1) = i(Z+1)^{-1}(Z-1) = i(Z-1)(Z+1)^{-1}  = W^{-1},
\]
hence $W \in U(2)$, and $\gamma$ maps $\BB M$ into $U(2)$.
Next, let $Z \in \BB T^{\pm}$, and let $W=-i(Z+1)(Z-1)^{-1}$. We compute:
\begin{multline*}
W^*W-1 = (Z^*-1)^{-1}(Z^*+1)(Z+1)(Z-1)^{-1} -1  \\
= (Z^*-1)^{-1} \bigl( (Z^*+1)(Z+1) - (Z^*-1)(Z-1) \bigr) (Z-1)^{-1}  \\
= 2 (Z^*-1)^{-1} (Z+ Z^*) (Z-1)^{-1}.
\end{multline*}
It follows that $W^*W-1$ is positive (or negative) definite if and only if
$Z+Z^*$ is. Thus $\gamma$ maps $\BB T^+$ and $\BB T^-$ into $\BB D^+$ and
$\BB D^-$ respectively.

Clearly, $\gamma$ on $\BB T^-$ has singularities along those elements that
have eigenvalue $1$, such as
$\bigl(\begin{smallmatrix} 1 & 0 \\ 0 & 1 \end{smallmatrix}\bigr) \in \BB T^-$.

Every element $Z \in \BB M \cup \BB T^+$ can be written as $Z=X+Y$ with
$X \in \BB M$ and $Y=0$ or $Y \in \HC$ such that $Y^*=Y$ and negative
definite. Each such self-adjoint negative definite $Y \in \HC$ is conjugate
to a diagonal matrix with negative real entries:
\[
Y = - UDU^*, \quad
D = \bigl(\begin{smallmatrix} \lambda_1 & 0 \\
  0 & \lambda_2 \end{smallmatrix}\bigr), \qquad
U \in U(2),\: \lambda_1,\lambda_2 >0.
\]
In either case, the determinant of $Z-1$ is
\[
N(Z-1) = N(X-UDU^*-1) = N(X'-D-1), \qquad \text{where } X'=U^*XU \in \BB M,
\]
and it is easy to verify that $N(X'-D-1) \ne 0$. This proves that $Z-1$
is invertible and $\gamma$ has no singularities on $\BB M \cup \BB T^+$.

Conversely, let $Z \in U(2)$, so that $Z^*=Z^{-1}$, and let $W=(Z-i)(Z+i)^{-1}$.
We compute:
\begin{multline*}
W^* = (Z^*-i)^{-1}(Z^*+i) = (Z^{-1}-i)^{-1}(Z^{-1}+i) = (1-iZ)^{-1}(1+iZ)  \\
= - (Z+i)^{-1}(Z-i) = - (Z-i)(Z+i)^{-1} = -W,
\end{multline*}
hence $W \in \BB M$, and $\gamma^{-1}$ maps $U(2)$ into $\BB M$.
Next, let $Z \in \BB D^{\pm}$, and let $W=(Z-i)(Z+i)^{-1}$. We compute:
\begin{multline*}
W+W^*= (Z-i)(Z+i)^{-1} + (Z^*-i)^{-1}(Z^*+i)  \\
= (Z^*-i)^{-1} \bigl( (Z^*-i)(Z-i) + (Z^*+i)(Z+i) \bigr) (Z+i)^{-1}  \\
= 2 (Z^*-i)^{-1} (Z^*Z-1) (Z+i)^{-1}.
\end{multline*}
It follows that $Z^*Z-1$ is positive (or negative) definite if and only if
$W+W^*$ is. Thus $\gamma^{-1}$ maps $\BB D^+$ and $\BB D^-$ into $\BB T^+$ and
$\BB T^-$ respectively.

The map $\gamma^{-1}$ on $U(2)$ has singularities along those elements that
have eigenvalue $-i$.
On the other hand, the pair $\{\lambda_1, \lambda_2\}$ in Proposition \ref{svd}
is uniquely determined by $Z \in \HC$.
It follows that $Z \in \BB D^+ \cup \BB D^-$ cannot have an eigenvalue of
absolute value 1, and hence $\gamma^{-1}$ has no singularities on
$\BB D^+ \cup \BB D^-$.
\end{proof}

We think of $\BB T^+$ and $\BB T^-$ as quaternionic analogues of
the upper and lower half planes.

Recall that the group $U(2,2)$ realized as a subgroup of $GL(2,\HC)$
as in \eqref{U(2,2)} acts on $U(2)$, $\BB D^+$ and $\BB D^-$ by fractional
linear transformations \eqref{conformal-action}.
Hence its conjugate by $\gamma$ acts on $\BB M$, $\BB T^+$ and $\BB T^-$.
Introduce
\begin{multline*}
U(2,2)' = \gamma^{-1} U(2,2) \gamma
= \Bigl\{ \bigl(\begin{smallmatrix} a & b \\
  c & d \end{smallmatrix}\bigr) \in GL(2,\HC) ;\:
\bigl(\begin{smallmatrix} a & b \\ c & d \end{smallmatrix}\bigr)
\bigl(\begin{smallmatrix} 0 & 1 \\ 1 & 0 \end{smallmatrix}\bigr)
\bigl(\begin{smallmatrix} a^* & c^* \\ b^* & d^* \end{smallmatrix}\bigr)
= \bigl(\begin{smallmatrix} 0 & 1 \\ 1 & 0 \end{smallmatrix}\bigr) \Bigr\}  \\
= \Biggl\{ \begin{pmatrix} a & b \\ c & d \end{pmatrix};\:
a,b,c,d \in \HC,\:
\begin{matrix} ab^*+ba^*=0 \\ cd^*+dc^*=0 \\ ad^*+bc^*=1 \end{matrix}
\Biggr\}.
\end{multline*}
(The relation between $U(2,2)$ and $U(2,2)'$ can be verified by observing
that $\gamma^*=\gamma^{-1}$ and
$\gamma \bigl(\begin{smallmatrix} 0 & 1 \\ 1 & 0 \end{smallmatrix}\bigr)
\gamma^* = \bigl(\begin{smallmatrix} 1 & 0 \\ 0 & -1 \end{smallmatrix}\bigr)$.)
We have a restatement of Proposition \ref{mobius-prop}:

\begin{prop}
The group $U(2,2)'$ acts on $\HC$ (with singularities)
by fractional linear transformations \eqref{conformal-action}
preserving $\BB M \subset \HC$ and tube domains $\BB T^+$ and $\BB T^-$.
\end{prop}

By definition, every element $Z$ of $\BB T^+$ can be written as
$Z=X+Y$ with $X \in \BB M$ and $Y \in \HC$ such that $Y^*=Y$ and negative
definite. Each such self-adjoint negative definite $Y \in \HC$ is conjugate
to a diagonal matrix with negative real entries:
\begin{equation}  \label{Y-diag}
Y = - UDU^*, \quad
D = \bigl(\begin{smallmatrix} \lambda_1 & 0 \\
  0 & \lambda_2 \end{smallmatrix}\bigr), \qquad
U \in U(2),\: \lambda_1,\lambda_2 >0.
\end{equation}
Denote by $D^{-\frac12}$ the matrix
$\Bigl(\begin{smallmatrix} \lambda_1^{-1/2} & 0 \\
  0 & \lambda_2^{-1/2} \end{smallmatrix}\Bigr)$
and define the {\em Poincar\'e metric} on $\BB T^+$ by declaring that,
for each tangent vector ${\bf v} \in T_Z\BB T^+ \simeq \HC$,
\[
|{\bf v}|'_Z = \tfrac12 \bigl\| D^{-\frac12} U^* {\bf v} U D^{-\frac12} \bigr\|_{HS}.
\]

\begin{prop}
The Poincar\'e metric on $\BB T^+$ is the pull-back of the Poincar\'e metric
on $\BB D^+$ via the Cayley transform $\gamma: \BB T^+ \to \BB D^+$.
The Poincar\'e metric on $\BB T^+$ is invariant under the action
of $U(2,2)'$ on $\BB T^+$ by fractional linear transformations
\eqref{conformal-action}.

Moreover, if a Riemannian metric on $\BB T^+$ is invariant under the action
of $U(2,2)'$ on $\BB T^+$ by fractional linear transformations
\eqref{conformal-action}, then this metric is a constant multiple of the
Poincar\'e metric.
\end{prop}

\begin{proof}
Let $|{\bf v}|''_Z$ denote the metric on $\BB T^+$ obtained by
pulling-back the Poincar\'e metric on $\BB D^+$ via the Cayley transform
$\gamma: \BB T^+ \to \BB D^+$.
By Proposition \ref{Poincere-invar-prop}, $|{\bf v}|''_Z$ is
$U(2,2)'$-invariant.
And by Proposition \ref{P-uniqueness-prop}, any other $U(2,2)'$-invariant
Riemannian metric on $\BB T^+$ is a constant multiple of $|{\bf v}|''_Z$.
It remains to show that $|{\bf v}|''_Z$ coincides with $|{\bf v}|'_Z$.

Observe that $\gamma$ maps
$-1 = \bigl(\begin{smallmatrix} -1 & 0 \\ 0 & -1 \end{smallmatrix}\bigr)
\in \BB T^+$ into $0 \in \BB D^+$, and the pull-back of the Poincar\'e metric
on $T_0 \BB D^+$ to $T_{-1} \BB T^+$ is $|{\bf v}|'_{-1}$.
Let $Y \in i\C^+ \subset \BB T^+$, then $Y \in \HC$ is self-adjoint and
negative definite, hence conjugate to a diagonal matrix with negative
real entries as in \eqref{Y-diag}.
Let $a = D^{-\frac12} U^* \in \HC$, then
$\bigl(\begin{smallmatrix} a & 0 \\ 0 & {a^*}^{-1} \end{smallmatrix}\bigr)
\in U(2,2)'$.
The fractional linear transformation on $\BB T^+$ corresponding to
$\bigl(\begin{smallmatrix} a & 0 \\ 0 & {a^*}^{-1} \end{smallmatrix}\bigr)$
maps $Y$ into $D^{-\frac12} U^* Y U D^{-\frac12} = -1$.
Hence $|{\bf v}|''_Z$ coincides with $|{\bf v}|'_Z$ along $iC^+$.
Finally, both $|{\bf v}|''_Z$ and $|{\bf v}|'_Z$ are invariant under
parallel translations of $\BB T^+$ by elements $X \in \BB M$ realized
by matrices
$\bigl(\begin{smallmatrix} 1 & X \\ 0 & 1 \end{smallmatrix}\bigr) \in U(2,2)'$.
Hence $|{\bf v}|''_Z$ coincides with $|{\bf v}|'_Z$ on all of $\BB T^+$.
\end{proof}

\begin{rem}
Similarly, one has an explicit description of a $U(2,2)'$-invariant
Riemannian metric on $\BB T^-$. (Such a metric is unique up to a constant
multiple.)
Every element $Z$ of $\BB T^-$ can be written as
$Z=X+Y$ with $X \in \BB M$ and $Y \in \HC$ such that $Y^*=Y$ and positive
definite. Each such self-adjoint positive definite $Y \in \HC$ is conjugate
to a diagonal matrix with positive real entries:
\[
Y = UDU^*, \quad
D = \bigl(\begin{smallmatrix} \lambda_1 & 0 \\
  0 & \lambda_2 \end{smallmatrix}\bigr), \qquad
U \in U(2),\: \lambda_1,\lambda_2 >0.
\]
Then the {\em Poincar\'e metric} on $\BB T^-$ can be defined by declaring that,
for each tangent vector ${\bf v} \in T_Z\BB T^- \simeq \HC$,
\[
|{\bf v}|'_Z = \tfrac12 \bigl\| D^{-\frac12} U^* {\bf v} U D^{-\frac12} \bigr\|_{HS}.
\]
\end{rem}

\section{The Spectral Distance and the Lipschitz Constant of Holomorphic Maps}
\label{spectral-sect}

Let $F: \BB D^+ \to \BB D^+$ be holomorphic. In this section we answer the following question: must $F$ be distance decreasing in the Poincar\'e metric? In other words, is it always true that, for any $P, Q \in \BB D^+$, $d(F(P),F(Q)) \le d(P,Q)$?

This question is answered in the negative, but $F$ is always Lipschitz with the sharp constant $\sqrt2$ (cf. Theorem \ref{sqrt2-thm} and Proposition \ref{sqrt2-sharp-prop} below). 

The mechanism is that $\BB D^+$ carries a second natural invariant metric,
a Finsler metric built out of the spectral norm \eqref{spectral-norm} rather
than the Hilbert--Schmidt norm, which {\em is} contracted by all holomorphic
self maps. It is this metric which is comparable to the Poincar\'e metric with the sharp
constant $\sqrt2$.

Recall that a Finsler metric is a norm on each tangent space, varying
smoothly with the base point, and its length distance is obtained by
integrating these norms along curves. The Kobayashi distance on a complex
manifold $X$ is the largest pseudodistance $k_X$ for which every holomorphic
map $\varphi:{\cal D}\to X$ is distance decreasing from the Poincar\'e
distance on ${\cal D}$, whereas the Carath\'eodory distance is
\[
c_X(P,Q)=\sup_{f:X\to{\cal D}\ \mathrm{holomorphic}}
 d_{\cal D}\bigl(f(P),f(Q)\bigr).
\]
For $\BB D^+$ these two intrinsic distances agree and will be identified
below with $\rho$; see \cite{H,Kob}.

\subsection{Another Invariant Finsler Metric}

We first record the Poincar\'e metric of Section \ref{Poincare-sect} in a form
that does not refer to a choice of singular value decomposition.
For $Z \in \BB D^+$ the matrices $1-ZZ^*$ and $1-Z^*Z$ are positive definite,
and we write $(1-ZZ^*)^{-\frac12}$, $(1-Z^*Z)^{-\frac12}$ for their positive
definite inverse square roots.

\begin{lem}  \label{intrinsic-metric-lem}
Let $Z \in \BB D^+$ be written as $Z = UDV^*$ as in Proposition \ref{svd},
and let ${\bf v} \in T_Z\BB D^+ \simeq \HC$. Then
\[
\bigl\| (1-D^2)^{-\frac12} (U^*{\bf v}V) (1-D^2)^{-\frac12} \bigr\|_{HS}
= \bigl\| (1-ZZ^*)^{-\frac12} {\bf v} (1-Z^*Z)^{-\frac12} \bigr\|_{HS},
\]
and the same identity holds with $\|\cdot\|_{HS}$ replaced by the spectral norm
$\|\cdot\|$ throughout.
\end{lem}

\begin{proof}
From $Z=UDV^*$ we get $ZZ^*=UD^2U^*$ and $Z^*Z = VD^2V^*$, hence
$(1-ZZ^*)^{-\frac12} = U(1-D^2)^{-\frac12}U^*$ and
$(1-Z^*Z)^{-\frac12} = V(1-D^2)^{-\frac12}V^*$. Therefore
\[
(1-ZZ^*)^{-\frac12} {\bf v} (1-Z^*Z)^{-\frac12}
= U \bigl[ (1-D^2)^{-\frac12} (U^*{\bf v}V) (1-D^2)^{-\frac12} \bigr] V^*,
\]
and both $\|\cdot\|_{HS}$ and $\|\cdot\|$ are unchanged by multiplication by
unitary matrices on either side.
\end{proof}

Thus, after transforming the tangent vector by the pair $(U,V)$, as in formula \eqref{|v|_Z-def}, the Poincar\'e metric is
\begin{equation}  \label{HS-Finsler}
|{\bf v}|_Z \;=\; \bigl\| (1-ZZ^*)^{-\frac12} {\bf v} (1-Z^*Z)^{-\frac12}
\bigr\|_{HS}.
\end{equation}
We now introduce its spectral companion.

\begin{df}  \label{spectral-Finsler-df}
The {\em spectral metric} on $\BB D^+$ is the Finsler metric
\begin{equation}  \label{sp-Finsler}
\|{\bf v}\|_Z \;=\; \bigl\| (1-ZZ^*)^{-\frac12} {\bf v} (1-Z^*Z)^{-\frac12}
\bigr\|,
\qquad {\bf v} \in T_Z \BB D^+ \simeq \HC ,
\end{equation}
where $\|\cdot\|$ is the spectral norm \eqref{spectral-norm}.
We denote by $\rho(P,Q)$ the associated length distance, i.e. the infimum of integrals
$\int_0^1 \|\gamma'(t)\|_{\gamma(t)}\,dt$ over piecewise ${\cal C}^1$ curves
$\gamma$ in $\BB D^+$ joining $P$ to $Q$.
\end{df}

Since $\|A\| \le \|A\|_{HS} \le \sqrt2\,\|A\|$ for $A \in \HC$, we have
\begin{equation}  \label{pointwise-comparison}
\|{\bf v}\|_Z \;\le\; |{\bf v}|_Z \;\le\; \sqrt2\,\|{\bf v}\|_Z ,
\qquad {\bf v} \in T_Z \BB D^+ .
\end{equation}
(The second inequality is the elementary bound $\lambda_1^2+\lambda_2^2 \le
2\max\{\lambda_1,\lambda_2\}^2$ for the singular values of a $2\times2$ matrix.)

\begin{prop}  \label{sp-invariance-prop}
The spectral metric \eqref{sp-Finsler} is invariant under the action of
$U(2,2)$ on $\BB D^+$ by fractional linear transformations
\eqref{conformal-action}. Consequently $\rho$ is a $U(2,2)$-invariant
pseudo-distance on $\BB D^+$.
\end{prop}

\begin{proof}
This is the proof of Proposition \ref{Poincere-invar-prop} verbatim.
The stabilizer \eqref{0-stabilizer} of $0$ acts on $T_0\BB D^+$ by
${\bf v} \mapsto a{\bf v}d^*$ with $a,d \in U(2)$, which preserves the spectral
norm; and for $g \in U(2,2)$ with $g(0)=UDV^*$ the tangent map
$\tilde g_*: T_0 \BB D^+ \to T_{g(0)}\BB D^+$ is
${\bf v} \mapsto U(1-D^2)^{\frac12}{\bf v}(1-D^2)^{\frac12}V^*$,
so that by Lemma \ref{intrinsic-metric-lem}
\[
\|\tilde g_* {\bf v}\|_{g(0)}
= \bigl\| (1-D^2)^{-\frac12}\bigl(U^* \tilde g_*{\bf v} V\bigr)(1-D^2)^{-\frac12}
\bigr\| = \|{\bf v}\| = \|{\bf v}\|_0 .
\]
The only property of the norm that was used is its invariance under
multiplication by unitary matrices.
\end{proof}

Theorem \ref{rho-formula-thm} below shows that $\rho$ is a distance rather than a pseudo-distance.

\subsection{Holomorphic Maps Contract the Spectral Metric}

The next proposition is the
infinitesimal form of the Schwarz Lemma of
Subsection \ref{Schwarz-subsect}; it is the reason the spectral norm, rather
than the Hilbert--Schmidt norm, is the correct one here.

\begin{prop}  \label{infinitesimal-SP-prop}
Let $F: \BB D^+ \to \BB D^+$ be holomorphic. Then for all $Z \in \BB D^+$ and
all ${\bf v} \in T_Z \BB D^+$,
\[
\| dF_Z({\bf v}) \|_{F(Z)} \;\le\; \|{\bf v}\|_Z .
\]
The same holds for holomorphic $f: \BB D^+ \to {\cal D}$, with the Poincar\'e
metric $|d\zeta|/(1-|\zeta|^2)$ on the unit disk ${\cal D}$.
\end{prop}

\begin{proof}
Fix $Z_1 \in \BB D^+$. By Lemma \ref{Mobius-D-lem} there are
$\psi, \tilde\psi \in U(2,2)$ with $\psi(0)=Z_1$ and $\tilde\psi(F(Z_1))=0$.
Then $G = \tilde\psi \circ F \circ \psi: \BB D^+ \to \BB D^+$ is holomorphic
with $G(0)=0$, so by Corollary \ref{Schwarz-cor0} (applied with $M=1$ and
${\cal V}=\HC$ equipped with the spectral norm) $\|G(W)\| \le \|W\|$ for all
$W \in \BB D^+$. Given ${\bf w} \in \HC$, ${\bf w} \ne 0$, apply this to
$W = t{\bf w}/\|{\bf w}\|$ with $t \in (0,1)$ real:
\[
\bigl\| G\bigl(t{\bf w}/\|{\bf w}\|\bigr) \bigr\| \le t .
\]
Since $G(t{\bf w}/\|{\bf w}\|) = t\,dG_0({\bf w})/\|{\bf w}\| + O(t^2)$,
dividing by $t$ and letting $t \to 0^+$ gives
$\|dG_0({\bf w})\| \le \|{\bf w}\|$, i.e. $\|dG_0({\bf w})\|_0 \le
\|{\bf w}\|_0$. As $\psi$ and $\tilde\psi$ are isometries of the spectral
metric by Proposition \ref{sp-invariance-prop}, the assertion follows with
${\bf v} = d\psi_0({\bf w})$.

For $f: \BB D^+ \to {\cal D}$ the same argument is used, replacing
$\tilde\psi$ by a fractional linear self-map of ${\cal D}$ carrying $f(Z_1)$
to $0$ and invoking Corollary \ref{Schwarz-cor1} instead.
\end{proof}

\begin{cor}  \label{rho-contracting-cor}
Every holomorphic $F: \BB D^+ \to \BB D^+$ satisfies
$\rho\bigl(F(P),F(Q)\bigr) \le \rho(P,Q)$ for all $P,Q \in \BB D^+$.
\end{cor}

\begin{proof}
$F$ maps a curve joining $P$ to $Q$ to a curve joining $F(P)$ to $F(Q)$ whose
$\rho$-length is, by Proposition \ref{infinitesimal-SP-prop}, at most the
$\rho$-length of the original curve.
\end{proof}


\begin{thm}  \label{rho-formula-thm}
For every $Z \in \BB D^+$,
\[
\rho(0,Z) \;=\; \frac12 \log \frac{1+\|Z\|}{1-\|Z\|}
\;=\; \operatorname{arctanh}\|Z\| .
\]
In particular $\rho$ is a genuine $U(2,2)$-invariant distance on $\BB D^+$.
\end{thm}

\begin{proof}
Write $Z = UDV^*$ as in Proposition \ref{svd}, with
$D = \bigl(\begin{smallmatrix} \lambda_1 & 0 \\ 0 & \lambda_2
\end{smallmatrix}\bigr)$ and $\lambda_1 = \|Z\| \ge \lambda_2 \ge 0$, and put
$t_j = \operatorname{arctanh}\lambda_j$.

{\em Upper bound.} Consider the curve
\[
\gamma(s) = U \begin{pmatrix} \tanh(st_1) & 0 \\ 0 & \tanh(st_2)\end{pmatrix}
V^*, \qquad s \in [0,1],
\]
which joins $0$ to $Z$ inside $\BB D^+$. Since
$\gamma(s) = U D(s) V^*$ with $D(s)$ diagonal, Lemma
\ref{intrinsic-metric-lem} gives
\[
\|\gamma'(s)\|_{\gamma(s)}
= \left\| \begin{pmatrix}
\dfrac{t_1\operatorname{sech}^2(st_1)}{1-\tanh^2(st_1)} & 0 \\[2mm]
0 & \dfrac{t_2 \operatorname{sech}^2(st_2)}{1-\tanh^2(st_2)}
\end{pmatrix}\right\|
= \left\| \begin{pmatrix} t_1 & 0 \\ 0 & t_2\end{pmatrix}\right\|
= t_1 .
\]
Hence $\rho(0,Z) \le \operatorname{Length}(\gamma) = t_1$.

{\em Lower bound.} Let $u \in \BB C^2$ and $w \in \BB C^2$ be unit vectors with
$Zw = \lambda_1 u$ (the leading left and right singular vectors of $Z$), and
define the linear map
\[
f: \HC \to \BB C, \qquad f(W) = u^* W w .
\]
For $W \in \BB D^+$ we have $|f(W)| \le \|W\| < 1$, so $f$ restricts to a
holomorphic map $\BB D^+ \to {\cal D}$; moreover $f(0)=0$ and
$f(Z) = \lambda_1$. By Proposition \ref{infinitesimal-SP-prop}, $f$ does not
increase lengths from the spectral metric to the Poincar\'e metric of
${\cal D}$. Therefore, for any curve $\gamma$ in $\BB D^+$ joining $0$ to $Z$,
\[
\operatorname{Length}_\rho(\gamma) \;\ge\;
\operatorname{Length}_{{\cal D}}(f \circ \gamma)
\;\ge\; d_{{\cal D}}\bigl(0,\lambda_1\bigr)
= \operatorname{arctanh}\lambda_1 = t_1 .
\]
Taking the infimum over $\gamma$ gives $\rho(0,Z) \ge t_1$.

Finally $\rho$ is a $U(2,2)$-invariant length pseudo-distance by Proposition
\ref{sp-invariance-prop}, and it separates points because $\rho(0,Z)>0$ for
$Z \ne 0$; by transitivity of the $U(2,2)$ action (Corollary \ref{hom-space-cor}) it separates all pairs of points.
\end{proof}

\begin{rem}  \label{rho-Kob-rem}
The infinitesimal Kobayashi norm on a complex manifold $X$ is
\[
k_X(Z;{\bf v})=\inf\left\{r^{-1};\ r>0,\ \exists\,\varphi:{\cal D}\to X
\text{ holomorphic with }\varphi(0)=Z,\ \varphi'(0)=r{\bf v}\right\}.
\]
This is the infinitesimal counterpart of the Kobayashi distance and is also
contracted by every holomorphic map. Following \cite{H,Kob},
Corollary \ref{rho-contracting-cor} and Theorem \ref{rho-formula-thm} identify
$\rho$ as the Kobayashi (equivalently, by the proof of the lower bound above,
the Carath\'eodory) distance of $\BB D^+$. The proof given
here is self-contained and uses only the Schwarz Lemma of Subsection
\ref{Schwarz-subsect}.
\end{rem}

\subsection{The Sharp Lipschitz Constant}

\begin{thm}  \label{sqrt2-thm}
Let $d$ denote the Poincar\'e distance of Section \ref{Poincare-sect}.
Then
\begin{equation}  \label{rho-d-comparison}
\rho(P,Q) \;\le\; d(P,Q) \;\le\; \sqrt2\,\rho(P,Q),
\qquad P,Q \in \BB D^+,
\end{equation}
and both inequalities are attained. Consequently, every holomorphic map
$F: \BB D^+ \to \BB D^+$ satisfies
\begin{equation}  \label{sqrt2-Lip}
d\bigl(F(P),F(Q)\bigr) \;\le\; \sqrt2\; d(P,Q),
\qquad P,Q \in \BB D^+ .
\end{equation}
\end{thm}

\begin{proof}
Both $d$ and $\rho$ are $U(2,2)$-invariant (Proposition
\ref{Poincere-invar-prop} and Proposition \ref{sp-invariance-prop}) and
$U(2,2)$ acts transitively on $\BB D^+$, so it suffices to prove
\eqref{rho-d-comparison} for $P=0$, $Q=Z$. Writing $Z=UDV^*$ as above and
$t_j=\operatorname{arctanh}\lambda_j$ with $t_1 \ge t_2 \ge 0$, Proposition
\ref{Poincare-dist-prop} and Theorem \ref{rho-formula-thm} give
\[
d(0,Z) = \sqrt{t_1^2+t_2^2}, \qquad \rho(0,Z)=t_1 ,
\]
and $t_1 \le \sqrt{t_1^2+t_2^2} \le \sqrt2\,t_1$. Equality holds on the left
exactly when $t_2=0$, i.e. when $Z$ is singular, and on the right exactly when
$t_1=t_2$, i.e. when $Z=rU$ with $r\in[0,1)$ and $U \in U(2)$.

For \eqref{sqrt2-Lip} combine \eqref{rho-d-comparison} with Corollary
\ref{rho-contracting-cor}:
\[
d\bigl(F(P),F(Q)\bigr) \le \sqrt2\,\rho\bigl(F(P),F(Q)\bigr)
\le \sqrt2\,\rho(P,Q) \le \sqrt2\, d(P,Q). \qedhere
\]
\end{proof}

\begin{prop}  \label{sqrt2-sharp-prop}
The constant $\sqrt2$ in \eqref{sqrt2-Lip} is optimal. Indeed, let
\[
F: \BB D^+ \to \BB D^+, \qquad
F \bigl(\begin{smallmatrix} z_{11} & z_{12} \\ z_{21} & z_{22}
\end{smallmatrix}\bigr) = z_{11} \cdot
\bigl(\begin{smallmatrix} 1 & 0 \\ 0 & 1 \end{smallmatrix}\bigr).
\]
Then $F$ is holomorphic, $F(0)=0$, and for every $\lambda \in [0,1)$, taking
$Z_\lambda = \bigl(\begin{smallmatrix} \lambda & 0 \\ 0 & 0
\end{smallmatrix}\bigr)$,
\[
d\bigl(F(0),F(Z_\lambda)\bigr) = \sqrt2 \; d(0,Z_\lambda).
\]
In particular a holomorphic map $F: \BB D^+ \to \BB D^+$ need not be distance
decreasing for the Poincar\'e metric.
\end{prop}

\begin{proof}
$F$ is linear, hence holomorphic, and
$\|F(Z)\| = |z_{11}| = |e_1^* Z e_1| \le \|Z\| < 1$ for $Z \in \BB D^+$, where
$e_1 = \bigl(\begin{smallmatrix} 1 \\ 0\end{smallmatrix}\bigr)$; so
$F(\BB D^+)\subset \BB D^+$. The singular values of $Z_\lambda$ are
$\lambda, 0$ and those of $F(Z_\lambda)= \lambda \cdot 1$ are $\lambda,\lambda$.
By Proposition \ref{Poincare-dist-prop},
\[
d(0,Z_\lambda) = \operatorname{arctanh}\lambda, \qquad
d\bigl(0, F(Z_\lambda)\bigr) = \sqrt2 \operatorname{arctanh}\lambda . \qedhere
\]
\end{proof}

\subsection{Two Rigidity Consequences}

The spectral metric also gives a clean analogue of the equality case in the
classical Schwarz Lemma.  We only need equality infinitesimally at the
origin, but in every tangent direction.

\begin{thm}  \label{infinitesimal-rigidity-thm}
Let $F:\BB D^+\to\BB D^+$ be holomorphic with $F(0)=0$.  Suppose that
\begin{equation}  \label{infinitesimal-rigidity-eqn}
\|dF_0({\bf v})\|=\|{\bf v}\| \qquad\hbox{for every }{\bf v}\in \HC,
\end{equation}
where $\|\cdot\|$ is the spectral norm.  Then there are $U,V\in U(2)$ such
that either
\[
F(Z)=UZV \qquad\hbox{for all }Z\in\BB D^+,
\]
or
\[
F(Z)=UZ^tV \qquad\hbox{for all }Z\in\BB D^+.
\]
In particular, $F$ is an automorphism of $\BB D^+$.
\end{thm}

\begin{proof}
By Proposition \ref{infinitesimal-SP-prop}, the spectral norm is precisely the
infinitesimal Kobayashi norm at the origin.  Put $T=dF_0$.  Condition
\eqref{infinitesimal-rigidity-eqn} says that $T$ is a complex-linear isometry
of the four-dimensional normed space $\HC\cong M_2(\BB C)$, and hence it is onto.  We use
the standard classification of surjective linear isometries of a full matrix
algebra: in the present $2\times2$ case it says simply that
\[
T(Z)=UZV \quad\hbox{or}\quad T(Z)=UZ^tV,
\]
with $U,V$ unitary (this is the special case of Kadison's Isometry Theorem \cite[Theorem 7]{Kad} that we need here).
Each of these maps is a linear automorphism
of the matrix ball.  Compose $F$ with the inverse of $T$.  We obtain a holomorphic
self-map $G$ of $\BB D^+$ with $G(0)=0$ and $dG_0=1$.  Cartan's uniqueness
theorem states that a holomorphic self-map of a bounded domain fixing a point and
having identity derivative there is the identity. This gives $G=1$ and proves the result.
\end{proof}

There is a complementary rigidity statement at the Euclidean boundary.  A
point $P\in\partial\BB D^+$ is smooth precisely when the top singular value
$1$ is simple; thus, if the singular values of $P$ are $\lambda_1\ge\lambda_2$, the
smooth case is $\lambda_2<1$.  Notice that this excludes the Shilov boundary
$U(2)$, where both singular values are equal to $1$.

\begin{thm}[Burns--Krantz rigidity]  \label{Burns-Krantz-thm}
Let $P\in\partial\BB D^+$ be a smooth boundary point and let
$F:\BB D^+\to\BB D^+$ be holomorphic.  If
\begin{equation}  \label{Burns-Krantz-eqn}
F(Z)=Z+o\bigl(\|Z-P\|_{HS}^{3}\bigr) \qquad\hbox{as }Z\to P
\end{equation}
from inside $\BB D^+$, then $F$ is the identity map.
\end{thm}

\begin{proof}
This is the Burns--Krantz theorem for bounded symmetric domains proved by
Ng and Rong \cite{NR}.  Our matrix ball is an irreducible bounded symmetric
domain in its Harish--Chandra realization, so their theorem applies at every
smooth boundary point.  The choice of the Hilbert--Schmidt norm in
\eqref{Burns-Krantz-eqn} is immaterial, since all norms on $M_2(\BB C)$ are
equivalent.
\end{proof}

\begin{rem}  \label{Shilov-BK-rem}
Theorem \ref{Burns-Krantz-thm} is deliberately not a statement at the Shilov
boundary.  Points of $U(2)$ have both singular values equal to $1$ and lie in
the nonsmooth part of $\partial\BB D^+$.  A one-point boundary Schwarz lemma
at such maximal-rank points is therefore a different question; the
barycenter construction of Section \ref{bary-sect} gives boundary information
of another kind.
\end{rem}

\begin{rem}  \label{Un-n-rem}
Everything in this section extends verbatim to
\[
\BB D^+_{(n)} = \{ Z \in \BB C^{n \times n} ;\: \|Z\|<1 \}
= U(n,n)/\bigl(U(n)\times U(n)\bigr)
\]
with the Poincar\'e metric $|{\bf v}|_Z = \|(1-ZZ^*)^{-1/2}{\bf v}
(1-Z^*Z)^{-1/2}\|_{HS}$ and the spectral metric \eqref{sp-Finsler}.
One obtains $\rho(0,Z)=\operatorname{arctanh}\|Z\|$,
$d(0,Z)=\bigl(\sum_{j=1}^n (\operatorname{arctanh}\lambda_j)^2\bigr)^{1/2}$,
$\rho \le d \le \sqrt n\,\rho$, and consequently
\[
d\bigl(F(P),F(Q)\bigr) \le \sqrt n \; d(P,Q)
\]
for every holomorphic $F: \BB D^+_{(n)} \to \BB D^+_{(n)}$; the constant
$\sqrt n$ is attained by $F(Z) = z_{11}\cdot 1$ at the pairs
$\bigl(0, \operatorname{diag}(\lambda,0,\dots,0)\bigr)$.
Thus the Poincar\'e metric is distance decreasing for all holomorphic self maps
precisely in the rank one case $n=1$, that is, precisely for the unit disk.
\end{rem}

\section{The Barycenter Method on $\BB D^+$}  \label{bary-sect}

The Schwarz Lemma of Subsection \ref{Schwarz-subsect} and the Schwarz--Pick
Theorems of Subsection \ref{Schwarz-Pick-subsect} all start from a
normalisation: one moves a distinguished point to the origin by a fractional
linear transformation. In one complex variable there is a second, purely
boundary-theoretic way of producing such a normalisation, due to Douady and
Earle \cite{DE}: a probability measure $\mu$ on the circle $\partial{\cal D}$
that is not too concentrated has a unique {\em conformal barycenter} $z = \operatorname{bar}(\mu)$ in
${\cal D}$, characterised by
\[
\int_{\partial {\cal D}} \frac{u - z}{1-\bar z u}\,d\mu(u) =0 ,
\]
and the assignment $\mu \mapsto \operatorname{bar}(\mu)$ commutes with the
action of $SU(1,1)$. Here ``not too concentrated'' means that no atom has
mass at least $\frac12$. If $h:\partial{\cal D}\to\partial{\cal D}$ is a
homeomorphism and $\nu_z$ is harmonic measure at $z$, Douady and Earle define
$Eh(z)=\operatorname{bar}(h_*\nu_z)$; this extension is conformally natural.
The same device -- in the Riemannian guise of a centre of mass for Busemann
functions -- is the engine of the barycenter method of Besson, Courtois and
Gallot \cite{BCG}.

In this section we develop this machinery on $\BB D^+$, with the Shilov
boundary $U(2)$ of Proposition \ref{Shilov-prop} playing the role of
$\partial {\cal D}$. Everything is self contained; the only inputs are the
Schwarz Lemma of Subsection \ref{Schwarz-subsect}, the M\"obius maps of Lemma
\ref{Mobius-D-lem} and the Poincar\'e metric of Section \ref{Poincare-sect}.
Throughout, $dU$ denotes the normalised Haar measure on the compact group
$U(2)$, so that $\int_{U(2)}dU=1$.

\subsection{A Cocycle Identity and Harmonic Measure on the Shilov Boundary}

Following the proof of Proposition \ref{mobius-prop} we write the action
\eqref{conformal-action} of $g = \bigl(\begin{smallmatrix} a & b \\
c & d \end{smallmatrix}\bigr) \in U(2,2)$ on $\HC$ as
\begin{equation}  \label{g-action}
g\cdot Z = K(g,Z)^{-1}\bigl(c^*+Zd^*\bigr), \qquad
K(g,Z) := a^*+Zb^* .
\end{equation}

\begin{lem}  \label{cocycle-lem}
Let $g = \bigl(\begin{smallmatrix} a & b \\ c & d \end{smallmatrix}\bigr)
\in U(2,2)$. Then for all $Z, W \in \overline{\BB D^+}$ the matrices $K(g,Z)$
and $K(g,W)$ are invertible and
\begin{equation}  \label{cocycle-eqn}
1 - (g\cdot Z)(g\cdot W)^* \;=\; K(g,Z)^{-1}\,\bigl(1-ZW^*\bigr)\,
K(g,W)^{-*},
\end{equation}
where $K^{-*}=(K^{*})^{-1}$. In particular $g$ restricts to a real analytic
diffeomorphism of $\overline{\BB D^+}$ preserving $\BB D^+$ and $U(2)$.
\end{lem}

\begin{proof}
Using $g\cdot W = K(g,W)^{-1}(c^*+Wd^*)$ we get
$(g\cdot W)^* = (c+dW^*)K(g,W)^{-*}$, hence
\[
1-(g\cdot Z)(g\cdot W)^*
= K(g,Z)^{-1}\Bigl[ K(g,Z)K(g,W)^* - (c^*+Zd^*)(c+dW^*)\Bigr]K(g,W)^{-*} .
\]
Expanding the bracket,
\begin{multline*}
(a^*+Zb^*)(a+bW^*) - (c^*+Zd^*)(c+dW^*)  \\
= (a^*a-c^*c) + (a^*b-c^*d)W^* + Z(b^*a-d^*c) + Z(b^*b-d^*d)W^* ,
\end{multline*}
and the three relations in \eqref{U(2,2)}, namely $a^*a=1+c^*c$, $d^*d=1+b^*b$
and $a^*b=c^*d$ (whence also $b^*a=d^*c$), turn this into $1-ZW^*$. This
proves \eqref{cocycle-eqn} once invertibility of $K$ is known.

For invertibility, let $Z \in \overline{\BB D^+}$ and suppose
$K(g,Z)^*x = (a+bZ^*)x = 0$ for some $x \in \BB C^2$. Taking $W=Z$ in the
computation above (which did not use invertibility) gives the identity
\[
K(g,Z)K(g,Z)^* = (1-ZZ^*) + (c^*+Zd^*)(c^*+Zd^*)^* ,
\]
so that
\[
0 = \|K(g,Z)^*x\|^2 = x^*(1-ZZ^*)x + \|(c+dZ^*)x\|^2 .
\]
Both terms are nonnegative because $\|Z\|\le1$, hence both vanish; in
particular $(c+dZ^*)x=0$. Together with $(a+bZ^*)x=0$ this says
$g\bigl(\begin{smallmatrix} x \\ Z^*x\end{smallmatrix}\bigr)=0$, and since $g$
is invertible we get $x=0$. Thus $K(g,Z)$ is invertible.
The last assertion of the lemma now follows from \eqref{cocycle-eqn} with $W=Z$, which shows
that $1-(g\cdot Z)(g\cdot Z)^*$ is positive definite, positive semidefinite of
rank $0$, respectively, exactly when $1-ZZ^*$ is; and $g^{-1} \in U(2,2)$
provides the inverse map.
\end{proof}

\begin{cor}  \label{det-cocycle-cor}
For $g \in U(2,2)$ and $Z,W \in \overline{\BB D^+}$,
\begin{equation}  \label{det-cocycle}
\det\bigl(1-(g\cdot Z)(g\cdot W)^*\bigr)
= \frac{\det(1-ZW^*)}{\det K(g,Z)\cdot \overline{\det K(g,W)}} .
\end{equation}
\end{cor}

Lemmas 60, 61 from \cite{FL1} imply:

\begin{lem}  \label{jacobian-lem}
Let $g \in U(2,2)$. The Jacobian of the diffeomorphism $U \mapsto g\cdot U$ of
$U(2)$ with respect to the Haar measure equals
$\bigl|\det K(g,U)\bigr|^{-4}$.
\end{lem}


\begin{df}  \label{harmonic-measure-df}
For $Z \in \BB D^+$ choose $g_Z \in U(2,2)$ with $g_Z\cdot 0 = Z$
(Lemma \ref{Mobius-D-lem}) and set
\[
\nu_Z := (g_Z)_* \bigl(dU\bigr),
\]
the push-forward of the normalised Haar measure of $U(2)$. We call $\nu_Z$ the
{\em harmonic measure} of $\BB D^+$ at $Z$.
\end{df}

\begin{prop}  \label{harmonic-measure-prop}
The measure $\nu_Z$ does not depend on the choice of $g_Z$.
Thus, $\nu_Z$ is a well defined Borel probability measure on $U(2)$; it satisfies
\begin{equation}  \label{Poisson-Szego}
d\nu_Z(U) = P(Z,U)\,dU, \qquad
P(Z,U) = \frac{\det(1-ZZ^*)^2}{\bigl|\det(1-ZU^*)\bigr|^{4}} ,
\end{equation}
and it is equivariant:
\begin{equation}  \label{nu-equivariance}
g_*\nu_Z = \nu_{g\cdot Z} \qquad \text{for all } g \in U(2,2).
\end{equation}
Moreover $\nu_0 = dU$.
\end{prop}

\begin{proof}
If $g_Z\cdot 0 = g_Z'\cdot 0 = Z$ then $g_Z^{-1}g_Z'$ lies in the stabiliser
\eqref{0-stabilizer} of $0$, and therefore acts on $U(2)$ by
$U \mapsto aUd^{-1}$ with $a,d \in U(2)$; such a map preserves Haar measure.
Hence $\nu_Z$ does not depend on the choice of $g_Z$, and
\eqref{nu-equivariance} holds because $g\,g_Z$ carries $0$ to $g\cdot Z$.
Clearly $\nu_0=dU$.

For \eqref{Poisson-Szego}, put $h=g_Z^{-1}$, so $h\cdot Z = 0$. By Lemma
\ref{jacobian-lem} the change of variables $U' = h\cdot U$ gives
$d\nu_Z(U') = |\det K(h,U')|^{-4}dU'$. Now \eqref{det-cocycle} with $W=Z$
yields $|\det K(h,Z)|^{2} = \det(1-ZZ^*)$, because
$\det(1-(h\cdot Z)(h\cdot Z)^*) = \det(1-0)=1$; and \eqref{det-cocycle} with
the pair $(Z,U')$ yields
\[
1 = \det\bigl(1-(h\cdot Z)(h\cdot U')^*\bigr)
= \frac{\det(1-Z U'^*)}{\det K(h,Z)\cdot\overline{\det K(h,U')}},
\]
so $\overline{\det K(h,U')} = \det(1-ZU'^*)/\det K(h,Z)$. Combining,
\[
\bigl|\det K(h,U')\bigr|^{-4}
= \frac{\bigl|\det K(h,Z)\bigr|^{4}}{\bigl|\det(1-ZU'^*)\bigr|^{4}}
= \frac{\det(1-ZZ^*)^2}{\bigl|\det(1-ZU'^*)\bigr|^{4}} = P(Z,U'). \qedhere
\]
\end{proof}

\begin{rem}  \label{Poisson-remarks}
(i) In particular $\int_{U(2)}P(Z,U)\,dU=1$, and $P$ is the Poisson-Szeg\H o
kernel of $\BB D^+$; cf. \cite{Hua,Ko}. The proof above derives it from the
cocycle \eqref{cocycle-eqn} alone.

(ii) In the ``half-plane'' picture of Section \ref{half-plane-sect} the kernel
takes a particularly transparent form. For $Z \in \BB D^+$ and $U \in U(2)$ put
$Y = ZU^*$ and let $C = (1-Y)^{-1}(1+Y)$ be the Cayley transform of $Y$; then
\[
P(Z,U) = \det \Bigl( \tfrac12 \bigl(C+C^*\bigr)\Bigr)^{2} ,
\]
because a direct computation gives
$\tfrac12(C+C^*) = (1-Y)^{-1}(1-YY^*)(1-Y^*)^{-1}$ while
$\det(1-YY^*)=\det(1-ZZ^*)$ and $\det(1-Y)=\det(1-ZU^*)$. Thus $P$ is the
square of the determinant of the ``imaginary part'', exactly as for the
Poisson kernel $\Im w/|w-x|^2$ of the upper half-plane.

(iii) It is worth comparing $P$ with the quaternionic Poisson kernel of
Theorem \ref{Poisson-thm}. If $X_0 \in {\cal B} = \BB H \cap \BB D^+$ and
$X \in \partial{\cal B} = \BB H \cap U(2)$, then $X_0X_0^* = N(X_0)\cdot1$ and
$\det(1-X_0X^*) = N(X-X_0)$, so that
\[
P(X_0,X) = \left( \frac{1-N(X_0)}{N(X-X_0)}\right)^{4},
\qquad\text{whereas}\qquad
\frac{1-N(X_0)}{N(X-X_0)^2}
\]
is the kernel of Theorem \ref{Poisson-thm} with $R=1$. The two are genuinely
different objects, and so are their measures: $\partial{\cal B}$ is a
three dimensional submanifold of the four dimensional Shilov boundary $U(2)$,
and $\nu_{X_0}$ is not supported on it. Theorem \ref{Poisson-thm} reproduces
functions harmonic for $\square$ on ${\cal B}$, whereas $P(Z,U)$ reproduces
functions holomorphic on $\BB D^+$ from their Shilov boundary values.
\end{rem}

\subsection{Transversality Sets}

Recall from Subsection 2.1 that elements $Z \in \HC$ are identified with the
planes $\bigl(\begin{smallmatrix} Z \\ 1\end{smallmatrix}\bigr)$ in the
Grassmannian $\operatorname{Gr}_2(\HC)$, and that $GL(2,\HC)$ acts on
$\operatorname{Gr}_2(\HC)$ linearly. Two planes may or may not intersect
non-trivially, and this is a $GL(2,\HC)$-invariant condition. On the Shilov
boundary it reads as follows.

\begin{df}  \label{transversality-df}
For $W \in U(2)$ set
\[
{\cal Z}_W := \bigl\{ U \in U(2) ;\: N(U-W)=0 \bigr\}
= \bigl\{ U \in U(2) ;\: \det(1-U^*W)=0 \bigr\} .
\]
\end{df}

Thus $U \in {\cal Z}_W$ exactly when the unitary matrices $U$ and $W$ agree on
some complex line in $\BB C^2$, equivalently when the corresponding planes
$\bigl(\begin{smallmatrix} U \\ 1\end{smallmatrix}\bigr)$ and
$\bigl(\begin{smallmatrix} W \\ 1\end{smallmatrix}\bigr)$ in
$\operatorname{Gr}_2(\HC)$ are not transverse.

\begin{lem}  \label{transversality-lem}
\begin{enumerate}
\item[(a)] ${\cal Z}_W$ is a closed subset of $U(2)$ of Haar measure zero.
\item[(b)] $g\cdot{\cal Z}_W = {\cal Z}_{g\cdot W}$ for every $g \in U(2,2)$
and $W \in U(2)$.
\end{enumerate}
\end{lem}

\begin{proof}
(a) Closedness is clear. The map $U \mapsto V=U^*W$ is a Haar measure
preserving bijection of $U(2)$ carrying ${\cal Z}_W$ onto
$\{V \in U(2) ;\: \det(1-V)=0\}$, which is the zero set of the real analytic
function $V \mapsto |\det(1-V)|^2$ on the connected real analytic manifold
$U(2)$. This function is not identically zero (its value at $V=-1$ is $4$),
so its zero set has measure zero.

(b) By the discussion above, $U \in {\cal Z}_W$ if and only if the planes
$\bigl(\begin{smallmatrix} U \\ 1\end{smallmatrix}\bigr)$,
$\bigl(\begin{smallmatrix} W \\ 1\end{smallmatrix}\bigr)$ have a nonzero
intersection, and $g \in U(2,2)\subset GL(2,\HC)$ acts by an invertible linear
map of $\HC = \BB C^4$, hence preserves this relation. Alternatively, apply
\eqref{cocycle-eqn} with $Z=U$, $W=W$ and take determinants:
$\det\bigl(1-(g\cdot U)(g\cdot W)^*\bigr)$ vanishes exactly when
$\det(1-UW^*)$ does.
\end{proof}

\begin{df}  \label{admissible-df}
A Borel probability measure $\mu$ on $U(2)$ is called {\em admissible} if
\begin{equation}  \label{star}
\theta(\mu) \;:=\; \sup_{W \in U(2)} \mu\bigl({\cal Z}_W\bigr) \;<\;
\frac12 .
\end{equation}
\end{df}

\begin{lem}  \label{admissible-lem}
\begin{enumerate}
\item[(a)] The supremum in \eqref{star} is attained; in particular $\mu$ is
admissible as soon as $\mu({\cal Z}_W)<\frac12$ for every $W$.
\item[(b)] Admissibility is preserved by the action of $U(2,2)$ on measures,
and $\theta(g_*\mu)=\theta(\mu)$.
\item[(c)] Every measure absolutely continuous with respect to $dU$ is
admissible; in particular every harmonic measure $\nu_Z$ is admissible, with
$\theta(\nu_Z)=0$.
\end{enumerate}
\end{lem}

\begin{proof}
(a) We check that $W \mapsto \mu({\cal Z}_W)$ is upper semicontinuous. Let
$W_n \to W$ and let $O \supseteq {\cal Z}_W$ be open. If ${\cal Z}_{W_n}
\not\subseteq O$ for infinitely many $n$, pick $U_n \in {\cal Z}_{W_n}
\setminus O$; by compactness a subsequence converges to some
$U \in U(2)\setminus O$, and $N(U-W)=\lim N(U_n-W_n)=0$ gives
$U \in {\cal Z}_W \subseteq O$, a contradiction. Hence
${\cal Z}_{W_n}\subseteq O$ for all large $n$, so
$\limsup_n\mu({\cal Z}_{W_n}) \le \mu(O)$; taking the infimum over such $O$
and using outer regularity of $\mu$ gives
$\limsup_n\mu({\cal Z}_{W_n}) \le \mu({\cal Z}_W)$. An upper semicontinuous
function on the compact group $U(2)$ attains its supremum.
(b) is immediate from Lemma \ref{transversality-lem}(b).
(c) follows from Lemma \ref{transversality-lem}(a).
\end{proof}

\subsection{The Energy of a Measure}

\begin{df}  \label{Busemann-df}
For $U \in U(2)$ define the {\em Busemann function}
$B_U : \BB D^+ \to \BB R$ by
\begin{equation}  \label{Busemann-eqn}
B_U(Z) \;=\; \log \frac{\bigl|\det(1-ZU^*)\bigr|^2}{\det(1-ZZ^*)}
\;=\; -\tfrac12 \log P(Z,U) .
\end{equation}
For a Borel probability measure $\mu$ on $U(2)$ the {\em energy} of $\mu$ is
\[
{\cal E}_\mu : \BB D^+ \to \BB R, \qquad
{\cal E}_\mu(Z) = \int_{U(2)} B_U(Z)\,d\mu(U) .
\]
\end{df}

Since the eigenvalues of $ZU^*$ have modulus at most $\|Z\|$ we have
\[
(1-\|Z\|)^2 \le |\det(1-ZU^*)| \le 4;
\]
therefore $B_U(Z)$ and all of its partial derivatives in $Z$ are bounded uniformly in
$U \in U(2)$ on compact subsets of $\BB D^+$. Consequently ${\cal E}_\mu$ is
real analytic on $\BB D^+$ and may be differentiated under the integral sign.

\begin{lem}  \label{energy-cocycle-lem}
For every $g \in U(2,2)$ and every $U \in U(2)$, $Z \in \BB D^+$,
\[
B_{g\cdot U}(g\cdot Z) = B_U(Z) - 2\log\bigl|\det K(g,U)\bigr| .
\]
Consequently, for every probability measure $\mu$ on $U(2)$,
\begin{equation}  \label{energy-cocycle}
{\cal E}_\mu \circ g \;=\; {\cal E}_{(g^{-1})_*\mu} \;+\; c(g,\mu)
\end{equation}
with a constant $c(g,\mu) = -2\int_{U(2)}\log|\det K(g,U)|\,
d\bigl((g^{-1})_*\mu\bigr)(U) \in \BB R$ that does not depend on $Z$.
\end{lem}

\begin{proof}
By Corollary \ref{det-cocycle-cor},
$\det\bigl(1-(g\cdot Z)(g\cdot Z)^*\bigr) = \det(1-ZZ^*)\,
|\det K(g,Z)|^{-2}$ and
\[
\bigl|\det(1-(g\cdot Z)(g\cdot U)^*)\bigr|^{2}
= |\det(1-ZU^*)|^{2}\,|\det K(g,Z)|^{-2}|\det K(g,U)|^{-2}.
\]
Dividing gives the first identity. Integrating it against $\mu$ after the
substitution $U = g\cdot U'$ gives \eqref{energy-cocycle}; the constant is
finite because $U \mapsto |\det K(g,U)|$ is continuous and nowhere zero on the
compact set $U(2)$ (Lemma \ref{cocycle-lem}).
\end{proof}

\begin{lem}  \label{expansion-lem}
For $v \in \HC$ and $U \in U(2)$,
\[
B_U(tv) = -2t\,\re \tr(vU^*)
+ t^2\Bigl( \tr(vv^*) - \re \tr\bigl((vU^*)^2\bigr)\Bigr) + O(t^3)
\qquad (t \to 0).
\]
Consequently, for a probability measure $\mu$ on $U(2)$ and
$M(\mu) := \int_{U(2)} U\,d\mu(U) \in \HC$,
\begin{align}
\frac{d}{dt}\Big|_{t=0} {\cal E}_\mu(tv) &= -2\,\re \tr\bigl(v\,M(\mu)^*\bigr),
\label{first-var} \\
\frac{d^2}{dt^2}\Big|_{t=0} {\cal E}_\mu(tv) &=
2\Bigl( \|v\|_{HS}^2 - \re \tr \int_{U(2)} (vU^*)^2\,d\mu(U) \Bigr).
\label{second-var}
\end{align}
\end{lem}

\begin{proof}
Expand $\log\det(1-A) = \tr\log(1-A) = -\tr A - \tfrac12\tr(A^2)+O(\|A\|^3)$.
With $A = tvU^*$ this gives
\[
\log|\det(1-tvU^*)|^2 = 2\re\log\det(1-tvU^*)
= -2t\,\re\tr(vU^*) - t^2\re\tr((vU^*)^2)+O(t^3),
\]
and with $A = t^2vv^*$ it gives
$-\log\det(1-t^2vv^*) = t^2\tr(vv^*)+O(t^4)$. Add.
\end{proof}

\begin{lem}  \label{trace-ineq-lem}
For $A \in \HC$ we have $\re\tr(A^2) \le \|A\|_{HS}^2$, with equality if and
only if $A$ is Hermitian.
\end{lem}

\begin{proof}
Let $\lambda_1,\lambda_2$ be the eigenvalues of $A$. By Schur's inequality
$|\lambda_1|^2+|\lambda_2|^2 \le \|A\|_{HS}^2$, with equality if and only if
$A$ is normal (cf. \cite{HornJohnson}). Hence
\[
\re\tr(A^2)=\re(\lambda_1^2+\lambda_2^2) \le |\lambda_1|^2+|\lambda_2|^2
\le \|A\|_{HS}^2 .
\]
Equality throughout forces $A$ normal and $\lambda_j^2 \ge 0$, i.e.
$\lambda_j \in \BB R$; a normal matrix with real eigenvalues is Hermitian.
Conversely, if $A=A^*$ then $\tr(A^2)=\tr(AA^*)=\|A\|_{HS}^2$.
\end{proof}

Combining Lemma \ref{expansion-lem} and Lemma \ref{trace-ineq-lem} (applied to
$A = vU^*$, for which $\|A\|_{HS}=\|v\|_{HS}$ since $U$ is unitary) we see that
\eqref{second-var} is always $\ge 0$, and that it vanishes precisely when
$vU^*$ is Hermitian for $\mu$-almost every $U$. This motivates:

\begin{df}  \label{Sigma-df}
For $v \in \HC$, $v \ne 0$, put
$\Sigma_v := \{ U \in U(2) ;\: vU^* = U v^* \}$.
\end{df}

\begin{lem}  \label{Sigma-lem}
Write $v = U_1 \Lambda V_1^*$ as in Proposition \ref{svd}, with
$\Lambda = \bigl(\begin{smallmatrix} s_1 & 0 \\ 0 & s_2\end{smallmatrix}\bigr)$,
$s_1 \ge s_2 \ge 0$, $s_1>0$. Then $\Sigma_v = U_1 T_v V_1^*$, where
\[
T_v = \begin{cases}
\bigl\{ \bigl(\begin{smallmatrix} \epsilon_1 & 0 \\
0 & \epsilon_2\end{smallmatrix}\bigr) ;\: \epsilon_j \in \{\pm1\}\bigr\}
& \text{if } s_1>s_2>0 \quad (\text{four points}), \\[1mm]
\bigl\{ H \in U(2) ;\: H=H^* \bigr\}
& \text{if } s_1=s_2 \quad (\text{two points and a $2$-sphere}), \\[1mm]
\bigl\{ \bigl(\begin{smallmatrix} \epsilon & 0 \\
0 & \zeta \end{smallmatrix}\bigr) ;\: \epsilon \in\{\pm1\},\ |\zeta|=1\bigr\}
& \text{if } s_2=0 \quad (\text{two circles}).
\end{cases}
\]
In every case there exist $W_1, W_2 \in U(2)$ with
$\Sigma_v \subseteq {\cal Z}_{W_1}\cup{\cal Z}_{W_2}$; in particular $\Sigma_v$
is a closed set of Haar measure zero.
\end{lem}

\begin{proof}
Put $T = V_1^*U^*U_1$, so that $U = U_1T^*V_1^*$ and $T \in U(2)$. Then
$vU^* = U_1\Lambda T U_1^*$ is Hermitian if and only if $\Lambda T$ is
Hermitian. Writing $T=(t_{ij})$, the condition $(\Lambda T)^*=\Lambda T$ reads
\begin{equation}  \label{LT-herm}
s_1 t_{11}, s_2 t_{22} \in \BB R, \qquad s_1t_{12} = s_2 \bar t_{21} .
\end{equation}
If $s_1>s_2$, unitarity of $T$ gives $|t_{12}|=|t_{21}|$ (both equal
$\sqrt{1-|t_{11}|^2}$), so \eqref{LT-herm} forces $t_{12}=t_{21}=0$; then
$t_{11},t_{22}$ are real of modulus one, i.e. $\pm1$ if $s_2>0$; if $s_2=0$
the second condition in \eqref{LT-herm} reads $s_1t_{12}=0$, so $t_{12}=0$,
hence $t_{21}=0$ by unitarity, $t_{11}=\pm1$ and $t_{22}$ is an arbitrary
unimodular number. If $s_1=s_2$ then $\Lambda T = s_1 T$ is Hermitian if and
only if $T$ is a Hermitian unitary matrix. In all three cases $T_v$ is stable
under $T \mapsto T^*$, which gives the stated description of
$\Sigma_v=U_1T_vV_1^*$.

For the covering statement put $W_1 = U_1V_1^*$ and note that for
$U = U_1TV_1^*$ with $T \in T_v$ one has $U^*W_1 = V_1 T^* V_1^*$, so
$U \in {\cal Z}_{W_1}$ if and only if $\det(1-T)=0$, i.e. $T$ has the
eigenvalue $1$. In the first case this covers all $T$ except $T=-1$, and
$U_1(-1)V_1^* = -W_1 \in {\cal Z}_{-W_1}$; take $W_2=-W_1$. In the second case
the same is true, since a Hermitian unitary matrix has eigenvalues $\pm1$ and
lacks the eigenvalue $1$ only for $T=-1$; again $W_2=-W_1$. In the third case
the branch $\epsilon=1$ is covered by $W_1$; for the branch $\epsilon=-1$ take
$W_2 = U_1\bigl(\begin{smallmatrix} -1 & 0 \\
0 & 1\end{smallmatrix}\bigr)V_1^*$, and then
$U^*W_2 = V_1\bigl(\begin{smallmatrix} 1 & 0 \\
0 & \bar\zeta\end{smallmatrix}\bigr)V_1^*$ has the eigenvalue $1$.
Finally $\Sigma_v$ has measure zero by Lemma \ref{transversality-lem}(a).
\end{proof}

\begin{lem}  \label{flats-lem}
For $A,B \in U(2)$ define
\[
\Phi_{A,B}: \BB R^2 \to \BB D^+, \qquad
\Phi_{A,B}(u_1,u_2) = A \begin{pmatrix} \tanh u_1 & 0 \\
0 & \tanh u_2 \end{pmatrix} B^* .
\]
Then $\Phi_{A,B}$ is an isometry of the Euclidean plane onto its image, for
the Poincar\'e distance $d$ of Section \ref{Poincare-sect}. Consequently:
\begin{enumerate}
\item[(i)] for every $v \in T_0\BB D^+$, writing $v = A
\bigl(\begin{smallmatrix} c_1 & 0 \\ 0 & c_2\end{smallmatrix}\bigr)B^*$ for $c_1,c_2\ge0$ as in
Proposition \ref{svd}, the curve $\sigma_v(t) = \Phi_{A,B}(tc_1,tc_2)$ is a
geodesic with $\sigma_v(0)=0$, $\sigma_v'(0)=v$ and speed $\|v\|_{HS}$, and
every geodesic through $0$ arises this way;
\item[(ii)] any two points of $\BB D^+$ lie on a common geodesic;
\item[(iii)] $d(0,Z) \ge \operatorname{arctanh}\|Z\|$, so closed $d$-balls in
$\BB D^+$ are compact.
\end{enumerate}
\end{lem}

\begin{proof}
Proposition \ref{Poincare-dist-prop} says exactly that
$d\bigl(0,\Phi_{A,B}(u)\bigr) = \|u\|$: the singular values of
$\Phi_{A,B}(u)$ are $|\tanh u_1|,|\tanh u_2|$ and
$\operatorname{arctanh}|\tanh u_j| = |u_j|$.

Next, for $s \in \BB R^2$ put $D_s = \operatorname{diag}(\tanh s_1,\tanh s_2)$,
$R_s = (1-D_s^2)^{-\frac12}$ and
\[
\tau_s = \begin{pmatrix} R_s & D_sR_s \\ D_sR_s & R_s \end{pmatrix}
\quad\in GL(2,\HC).
\]
This is the matrix of Lemma \ref{Mobius-D-lem} with $U_0=V_0=1$ and $D_0=D_s$.
It lies in $U(2,2)$: the three relations in \eqref{U(2,2)} read
$a^*a=d^*d=R_s^2$, $c^*c=b^*b=R_s^2D_s^2$ and $a^*b=c^*d=R_s^2D_s$, and they
hold because $1+R_s^2D_s^2 = 1 + \frac{D_s^2}{1-D_s^2} = R_s^2$. Note that
only the reality and diagonality of $D_s$, together with $\|D_s\|<1$, were
used, so $s \in \BB R^2$ may have components of either sign.
The associated fractional linear transformation is
$W \mapsto R_s(W+D_s)(D_sW+1)^{-1}R_s^{-1}$, and if
$W = \operatorname{diag}(\tanh u_1,\tanh u_2)$ then every matrix occurring
here is diagonal, the factors $R_s^{\pm1}$ cancel, and the addition formula
for $\tanh$ gives
\[
\tau_s \cdot \Phi_{1,1}(u) = \Phi_{1,1}(u+s), \qquad u,s \in \BB R^2 .
\]
Finally, $k_{A,B} := \bigl(\begin{smallmatrix} A & 0 \\ 0 & B
\end{smallmatrix}\bigr)$ lies in the stabilizer \eqref{0-stabilizer} of $0$ in
$U(2,2)$ and acts by $W \mapsto AWB^*$, so that
$\Phi_{A,B} = k_{A,B}\cdot\Phi_{1,1}$ and
\[
\bigl(k_{A,B}\,\tau_s\,k_{A,B}^{-1}\bigr)\cdot \Phi_{A,B}(u)
= \Phi_{A,B}(u+s) .
\]
Since every element of $U(2,2)$ acts as an isometry of $d$ (Proposition
\ref{Poincere-invar-prop}), taking $s=-u'$ gives
\[
d\bigl(\Phi_{A,B}(u),\Phi_{A,B}(u')\bigr)
= d\bigl(\Phi_{A,B}(u-u'),\Phi_{A,B}(0)\bigr)
= d\bigl(\Phi_{A,B}(u-u'),0\bigr) = \|u-u'\| ,
\]
as claimed.

(i) $\sigma_v$ is the image under an isometry of a straight line, hence a
geodesic; it has $\sigma_v'(0)=v$ and speed $\sqrt{c_1^2+c_2^2}=\|v\|_{HS}$.
Since a geodesic is determined by its initial position and velocity, every
geodesic through $0$ is of this form.
(ii) Move the first point to $0$ by $U(2,2)$ (Corollary after Lemma
\ref{Mobius-D-lem}) and apply (i) to the second.
(iii) Immediate from Proposition \ref{Poincare-dist-prop}.
\end{proof}

\begin{prop}  \label{convexity-prop}
Let $\mu$ be a probability measure on $U(2)$ and let $\gamma$ be a non-constant
geodesic of the Poincar\'e metric in $\BB D^+$. Then $t \mapsto {\cal E}_\mu(\gamma(t))$ is
convex. If moreover $\mu$ is admissible, then it is strictly convex; more
precisely $\frac{d^2}{dt^2}{\cal E}_\mu(\gamma(t))>0$ for all $t$.
\end{prop}

\begin{proof}
By Lemma \ref{flats-lem}(i) the geodesics of the Poincar\'e metric through $0$
are exactly the curves
\begin{equation}  \label{geodesic-eqn}
\sigma(t) = A \begin{pmatrix} \tanh(c_1t) & 0 \\ 0 & \tanh(c_2t)\end{pmatrix}
B^*, \qquad A,B \in U(2),\ c_1,c_2 \ge 0 ,
\end{equation}
which have speed $\sqrt{c_1^2+c_2^2}$; as $t \mapsto \tanh(ct)$ has vanishing
second derivative at $t=0$, such a curve satisfies $\sigma''(0)=0$ as a curve
in the vector space $\HC$. Hence, for any smooth function $f$ on $\BB D^+$,
\begin{equation}  \label{Hess-at-0}
\frac{d^2}{dt^2}\Big|_{t=0} f\bigl(\sigma(t)\bigr)
= \frac{d^2}{dt^2}\Big|_{t=0} f\bigl(t\,\sigma'(0)\bigr) .
\end{equation}

Now let $\gamma$ be any geodesic and $t_0$ any parameter value. Choose
$g \in U(2,2)$ with $g\cdot 0 = \gamma(t_0)$. Since $g$ is an isometry of the
Poincar\'e metric (Proposition \ref{Poincere-invar-prop}), $\sigma := g^{-1}
\circ\gamma(\,\cdot + t_0)$ is a geodesic through $0$, hence of the form
\eqref{geodesic-eqn}. By \eqref{energy-cocycle} and \eqref{Hess-at-0},
\[
\frac{d^2}{dt^2}\Big|_{t=t_0}{\cal E}_\mu(\gamma(t))
= \frac{d^2}{dt^2}\Big|_{t=0}{\cal E}_{(g^{-1})_*\mu}\bigl(\sigma(t)\bigr)
= \frac{d^2}{dt^2}\Big|_{t=0}{\cal E}_{(g^{-1})_*\mu}\bigl(tv\bigr),
\qquad v = \sigma'(0)\ne 0 .
\]
By \eqref{second-var} and Lemma \ref{trace-ineq-lem} the last quantity equals
\[
2\int_{U(2)} \Bigl( \|v\|_{HS}^2 - \re\tr\bigl((vU^*)^2\bigr)\Bigr)
\,d\bigl((g^{-1})_*\mu\bigr)(U) \;\ge\; 0 ,
\]
which proves convexity. It vanishes if and only if the (nonnegative,
continuous) integrand vanishes $(g^{-1})_*\mu$-almost everywhere, i.e. if and
only if $(g^{-1})_*\mu(\Sigma_v)=1$. If $\mu$ is admissible then so is
$(g^{-1})_*\mu$ (Lemma \ref{admissible-lem}(b)), and by Lemma \ref{Sigma-lem}
\[
(g^{-1})_*\mu(\Sigma_v) \le (g^{-1})_*\mu({\cal Z}_{W_1})
+ (g^{-1})_*\mu({\cal Z}_{W_2}) \le 2\theta(\mu) <1 ,
\]
a contradiction. Hence the second derivative is strictly positive.
\end{proof}

\begin{prop}  \label{coercivity-prop}
Let $\mu$ be admissible and put $\delta_0 = 1-2\theta(\mu)>0$. Then there is
$t_0>0$ such that
\[
{\cal E}_\mu(Z) \;\ge\; {\cal E}_\mu(0) + \delta_0\, d(0,Z)
\qquad \text{whenever } d(0,Z) \ge t_0 .
\]
In particular ${\cal E}_\mu$ is bounded below, and its sublevel sets are
compact.
\end{prop}

\begin{proof}
{\em Step 1: a universal lower bound along geodesic rays.}
Let $\sigma$ be as in \eqref{geodesic-eqn} with $c_1,c_2\ge 0$ and $c_1^2+c_2^2=1$ (a unit speed
geodesic ray from $0$), and set $\Lambda(t) = \operatorname{diag}
(\tanh(c_1t),\tanh(c_2t))$ and, for $U \in U(2)$, $V = B^*U^*A \in U(2)$.
Since $\det\bigl(1-\sigma(t)\sigma(t)^*\bigr) = \prod_j
\operatorname{sech}^2(c_jt)$ and
$\det(1-\sigma(t)U^*)=\det(1-\Lambda(t)V)$,
\begin{equation}  \label{B-along-ray}
B_U(\sigma(t)) = 2\log\cosh(c_1t)+2\log\cosh(c_2t)
+ 2\log\bigl|\det(1-\Lambda(t)V)\bigr| .
\end{equation}
Let $z_1,z_2$ be the eigenvalues of $\Lambda(t)V$ and $\lambda_j=\tanh(c_jt)$.
Then $|z_j| \le \|\Lambda(t)V\| = \lambda_1$ and
$|z_1z_2| = |\det \Lambda(t)V| = \lambda_1\lambda_2$. We claim
\begin{equation}  \label{det-lower}
\bigl|\det(1-\Lambda(t)V)\bigr| = |1-z_1||1-z_2| \ge (1-\lambda_1)(1-\lambda_2).
\end{equation}
Indeed $|1-z_j| \ge 1-|z_j|$, and if $\lambda_2=0$ then one of the $z_j$
vanishes and \eqref{det-lower} reads $|1-z_1| \ge 1-\lambda_1$; if
$\lambda_2>0$ then $\lambda_2 \le |z_j| \le \lambda_1$ and
$(1-|z_1|)(1-|z_2|) = h(|z_1|)$ where
$h(x) = (1-x)(1-\lambda_1\lambda_2/x)$ is concave on $[\lambda_2,\lambda_1]$
and takes the value $(1-\lambda_1)(1-\lambda_2)$ at both endpoints.
Since $\cosh(ct)\bigl(1-\tanh(ct)\bigr) = \cosh(ct)-\sinh(ct) = e^{-ct}$, we
have $2\log\cosh(ct)+2\log(1-\tanh(ct)) = -2ct$, and \eqref{B-along-ray}
together with \eqref{det-lower} gives the universal bound
\begin{equation}  \label{universal-bound}
B_U(\sigma(t)) \;\ge\; -2(c_1+c_2)\,t \qquad
\text{for all } t\ge0 \text{ and all } U \in U(2).
\end{equation}

{\em Step 2: the asymptotic slope.} Suppose $U \notin {\cal Z}_{W_\sigma}$
where $W_\sigma := AB^*$. Since $1-V = B^*(B-U^*A)$ and, using $(AB^*-U)B = -(UB-A)$,
\begin{multline*}
\det(1-V) = \det(B^*)\det\bigl(U^*(UB-A)\bigr)  \\
= \det(B^*)\det(U^*)\det(AB^*-U)\det(B) = \det(U^*)\,\det(W_\sigma-U),
\end{multline*}
the condition $U\notin{\cal Z}_{W_\sigma}$ says exactly $\det(1-V)\ne0$. As $t \to\infty$ we have
$\Lambda(t) \to \Lambda_\infty := \operatorname{diag}(\chi_1,\chi_2)$ with
$\chi_j=1$ if $c_j>0$ and $\chi_j=0$ if $c_j=0$. If $c_1,c_2>0$ then
$\det(1-\Lambda(t)V) \to \det(1-V)\ne0$. If $c_2=0$ (so $c_1=1$) then
$\det(1-\Lambda(t)V) \to 1-v_{11}$, and $v_{11}=1$ would force
$V = \operatorname{diag}(1,\zeta)$ and hence $\det(1-V)=0$; so again the limit
is nonzero. In either case the last term of \eqref{B-along-ray} stays bounded
while the first two contribute $2(c_1+c_2)t+O(1)$, so
\begin{equation}  \label{slope}
\lim_{t\to\infty} \frac{B_U(\sigma(t))}{t} = 2(c_1+c_2),
\qquad U \notin {\cal Z}_{W_\sigma} .
\end{equation}

{\em Step 3: the slope of the energy.} The functions
$f_t(U) := B_U(\sigma(t))/t + 2(c_1+c_2)$ are nonnegative by
\eqref{universal-bound}, so Fatou's lemma together with \eqref{slope} gives
\begin{align*}
\liminf_{t\to\infty}\frac{{\cal E}_\mu(\sigma(t))}{t} + 2(c_1+c_2)
&\ge \int_{U(2)} \liminf_{t\to\infty} f_t \,d\mu
\;\ge\; 4(c_1+c_2)\,\mu\bigl(U(2)\setminus {\cal Z}_{W_\sigma}\bigr),
\end{align*}
whence, using $c_1+c_2 \ge \sqrt{c_1^2+c_2^2}=1$,
\begin{equation}  \label{slope-bound}
s(\sigma) := \lim_{t\to\infty}\frac{{\cal E}_\mu(\sigma(t))-{\cal E}_\mu(0)}{t}
\;\ge\; 2(c_1+c_2)\bigl(1-2\mu({\cal Z}_{W_\sigma})\bigr) \;\ge\; 2\delta_0 .
\end{equation}
(The limit exists in $(-\infty,+\infty]$ because
$q_t(\sigma):=\bigl({\cal E}_\mu(\sigma(t))-{\cal E}_\mu(0)\bigr)/t$ is
nondecreasing in $t$, ${\cal E}_\mu\circ\sigma$ being convex by Proposition
\ref{convexity-prop}.)

{\em Step 4: uniformity.} Let ${\cal S}$ be the (compact) set of unit vectors
of $\bigl(T_0\BB D^+, |\cdot|_0\bigr)$; by Lemma \ref{flats-lem}(i) each
$v \in {\cal S}$ determines a unit speed geodesic ray
$\sigma_v(t)=\exp_0(tv)$ from $0$, and $(v,t)\mapsto \sigma_v(t)$ is
continuous, being the Riemannian exponential map. Put
$m(t) = \min_{v \in {\cal S}} q_t(\sigma_v)$, which is attained by continuity
and is nondecreasing in $t$. We claim $\lim_{t\to\infty}m(t)\ge2\delta_0$.
Indeed, choose $t_n \to \infty$ with $m(t_n) \to L := \lim_t m(t)$ and
minimisers $v_n$ with $q_{t_n}(\sigma_{v_n}) = m(t_n)$; passing to a
subsequence, $v_n \to v_* \in {\cal S}$. For fixed $s>0$ and $t_n \ge s$,
monotonicity gives $q_s(\sigma_{v_n}) \le q_{t_n}(\sigma_{v_n}) = m(t_n)$;
letting $n\to\infty$ yields $q_s(\sigma_{v_*}) \le L$, and letting
$s \to \infty$ yields $2\delta_0 \le s(\sigma_{v_*}) \le L$ by
\eqref{slope-bound}.

Hence there is $t_0$ with $m(t)\ge\delta_0$ for $t \ge t_0$. Since every
$Z \in \BB D^+$ with $d(0,Z)=t\ge t_0$ is of the form $\sigma_v(t)$ for some
$v \in {\cal S}$ (Lemma \ref{flats-lem}(i)), we get
${\cal E}_\mu(Z)-{\cal E}_\mu(0) \ge \delta_0 t$, as claimed. Since closed
$d$-balls in $\BB D^+$ are compact (Lemma \ref{flats-lem}(iii)), the sublevel
sets of ${\cal E}_\mu$ are compact.
\end{proof}

\begin{thm}  \label{barycenter-thm}
Let $\mu$ be an admissible probability measure on $U(2)$. Then ${\cal E}_\mu$
attains its minimum at a unique point ${\cal B}(\mu) \in \BB D^+$, which is
also the unique critical point of ${\cal E}_\mu$. It is characterised by the
{\em barycenter equation}
\begin{equation}  \label{bary-eqn}
\int_{U(2)} (U-Z)\bigl(1-Z^*U\bigr)^{-1}\,d\mu(U) \;=\; 0,
\qquad Z = {\cal B}(\mu) .
\end{equation}
In particular ${\cal B}(\mu)=0$ if and only if $\int_{U(2)}U\,d\mu(U)=0$.
\end{thm}

\begin{proof}
Existence follows from Proposition \ref{coercivity-prop} (a continuous
function with compact sublevel sets attains its minimum). Any two points of
$\BB D^+$ are joined by a geodesic (Lemma \ref{flats-lem}(ii)), and
${\cal E}_\mu$ is strictly convex along every
geodesic by Proposition \ref{convexity-prop}; hence if $Z_1$ is a critical
point and $Z_2 \ne Z_1$, the function $t \mapsto {\cal E}_\mu(\gamma(t))$
along a geodesic $\gamma$ from $Z_1$ to $Z_2$ is strictly convex with
vanishing derivative at $t=0$, so ${\cal E}_\mu(Z_2)>{\cal E}_\mu(Z_1)$. Thus
there is at most one critical point, and it is the unique minimum.

For \eqref{bary-eqn}, let $Z \in \BB D^+$ and write $Z=U_0D_0V_0^*$ as in
Proposition \ref{svd}. By Lemma \ref{Mobius-D-lem} the element $g\in U(2,2)$
with
\[
g\cdot W = (1-D_0^2)^{-\frac12}U_0^*\,(W-Z)(1-Z^*W)^{-1}\,V_0(1-D_0^2)^{\frac12}
\]
carries $Z$ to $0$. By \eqref{energy-cocycle}, $Z$ is a critical point of
${\cal E}_\mu$ if and only if $0$ is a critical point of ${\cal E}_{g_*\mu}$,
which by \eqref{first-var} happens if and only if
$\int_{U(2)}U\,d(g_*\mu)(U)=\int_{U(2)}g\cdot U\,d\mu(U)=0$, i.e.
\[
(1-D_0^2)^{-\frac12}U_0^*\left[\int_{U(2)}(U-Z)(1-Z^*U)^{-1}d\mu(U)\right]
V_0(1-D_0^2)^{\frac12} = 0 .
\]
The outer factors are invertible and independent of $U$, so this is equivalent
to \eqref{bary-eqn}. (Note $1-Z^*U$ is invertible because
$\|Z^*U\|=\|Z\|<1$.)
\end{proof}

\begin{df}  \label{bary-df}
${\cal B}(\mu) \in \BB D^+$ is called the {\em barycenter} of the admissible
measure $\mu$.
\end{df}

The terminology reflects the two descriptions above: ${\cal B}(\mu)$ is the
unique minimizer of the averaged Busemann potential ${\cal E}_\mu$, and,
after moving it to the origin by an element of $U(2,2)$, the ordinary matrix
first moment of the transported measure vanishes. At the origin this is just
the familiar condition $M(\mu)=\int U\,d\mu(U)=0$.

\begin{rem}  \label{sharpness-rem}
Neither half of Theorem \ref{barycenter-thm} survives without \eqref{star},
and the constant $\frac12$ is the right one.

(a) {\em Coercivity fails.} Let $\mu = a\delta_{W} + (1-a)\mu'$ with
$W \in U(2)$ and $\mu'$ supported off ${\cal Z}_W$. Along the radial ray
$\sigma(t) = \tanh(t/\sqrt2)\,W$ (so $c_1=c_2=1/\sqrt2$) the computation in
Step 2 of Proposition \ref{coercivity-prop} gives
$\lim_{t\to\infty}B_U(\sigma(t))/t = 2(c_1+c_2)=2\sqrt2$ for $\mu'$-almost
every $U$, while for $U=W$ one has $V=1$, so that
$\det(1-\Lambda(t)V)=(1-\lambda_1)(1-\lambda_2)$ and, by the identity
$2\log\cosh(ct)+2\log(1-\tanh(ct))=-2ct$ used in Step 1,
$\lim_{t\to\infty}B_W(\sigma(t))/t = -2\sqrt2$. Therefore
\[
\lim_{t\to\infty}\frac{{\cal E}_\mu(\sigma(t))}{t} = 2\sqrt2\,(1-2a),
\]
which is negative as soon as $a>\frac12$: then ${\cal E}_\mu$ is unbounded
below and has no critical point at all. This is the exact analogue of the
Douady-Earle condition ``no atom of mass $\ge\frac12$'' \cite{DE}; note
however that in $\BB D^+$ it is the mass of the {\em three dimensional} set
${\cal Z}_W$, not merely of the point $W$, that has to be controlled.

(b) {\em Strict convexity fails.} Let $W \in U(2)$ and let $\mu$ be the
rotation invariant probability measure on the two-sphere
$\{HW ;\: H=H^*=H^{-1},\ \tr H = 0\}$, which by Lemma \ref{Sigma-lem} is
contained in $\Sigma_{W}$. Then $M(\mu)=0$, because $H \mapsto -H$ preserves
$\mu$, so $0$ satisfies the barycenter equation \eqref{bary-eqn}; but by
Lemma \ref{trace-ineq-lem} the second variation \eqref{second-var} vanishes at
$v=W$, so $0$ is a degenerate critical point. Here $\mu({\cal Z}_W)=1$, since
$U^*W = W^*HW$ has the eigenvalue $1$ for every such $U=HW$.
\end{rem}

\subsection{Naturality and the Barycentric Extension}

\begin{thm}  \label{naturality-thm}
For every admissible probability measure $\mu$ on $U(2)$ and every
$g \in U(2,2)$,
\[
{\cal B}\bigl(g_*\mu\bigr) = g \cdot {\cal B}(\mu) .
\]
Moreover ${\cal B}(\nu_Z)=Z$ for every $Z \in \BB D^+$.
\end{thm}

\begin{proof}
By \eqref{energy-cocycle}, ${\cal E}_{g_*\mu}\circ g = {\cal E}_\mu +
c(g,g_*\mu)$, so $g$ carries the minimum of ${\cal E}_\mu$ to the minimum of
${\cal E}_{g_*\mu}$; admissibility is preserved by Lemma
\ref{admissible-lem}(b). Since $\nu_0=dU$ is invariant under
$U \mapsto e^{i\vartheta}U$, we get $M(\nu_0)=e^{i\vartheta}M(\nu_0)$ for all
$\vartheta$, hence $M(\nu_0)=0$ and ${\cal B}(\nu_0)=0$ by Theorem
\ref{barycenter-thm}; the general case follows from \eqref{nu-equivariance}
and naturality.
\end{proof}

\begin{df}  \label{extension-df}
Let $\Psi: U(2)\to U(2)$ be a Borel map such that $\Psi_*\nu_Z$ is admissible
for every $Z \in \BB D^+$. The {\em barycentric extension} of $\Psi$ is
\[
E\Psi : \BB D^+ \to \BB D^+, \qquad
E\Psi(Z) := {\cal B}\bigl(\Psi_*\nu_Z\bigr) .
\]
\end{df}

The admissibility hypothesis holds, for instance, whenever $\Psi^{-1}$ carries
Haar null sets to Haar null sets: then $\Psi_*\nu_Z({\cal Z}_W) =
\nu_Z(\Psi^{-1}{\cal Z}_W)=0$ by Lemma \ref{transversality-lem}(a). It also
holds whenever $\sup_W \nu_Z(\Psi^{-1}{\cal Z}_W)<\frac12$ for all $Z$.
Here is a criterion that is easy to check in practice.

\begin{lem}  \label{analytic-admissible-lem}
Let $\Psi: U(2)\to U(2)$ be real analytic and suppose that for every
$W \in U(2)$ the function $U \mapsto \bigl|N\bigl(\Psi(U)-W\bigr)\bigr|^2$ is
not identically zero on $U(2)$. Then $\Psi_*\nu_Z$ is admissible for every
$Z\in\BB D^+$, with $\theta(\Psi_*\nu_Z)=0$.
\end{lem}

\begin{proof}
$\Psi^{-1}({\cal Z}_W)$ is the zero set of a real analytic function on the
connected real analytic manifold $U(2)$ that does not vanish identically,
hence has Haar measure zero; and $\nu_Z$ is absolutely continuous with respect
to Haar measure by \eqref{Poisson-Szego}.
\end{proof}

\begin{thm}  \label{extension-thm}
Let $\Psi$ be as in Definition \ref{extension-df}. Then:
\begin{enumerate}
\item[(a)] $E\Psi$ is real analytic on $\BB D^+$;
\item[(b)] {\em (conformal naturality)} for all $g,h \in U(2,2)$,
\[
E\bigl( g\circ\Psi\circ h^{-1}\bigr) \;=\; g \circ (E\Psi) \circ h^{-1} ,
\]
where $g,h$ act on $U(2)$ and on $\BB D^+$ by \eqref{conformal-action};
\item[(c)] $E(\operatorname{id})=\operatorname{id}$, and
$E\bigl(g|_{U(2)}\bigr)= g|_{\BB D^+}$ for every $g \in U(2,2)$.
\end{enumerate}
\end{thm}

\begin{proof}
(b) By \eqref{nu-equivariance}, $\nu_{h\cdot Z} = h_*\nu_Z$, so
$(g\circ\Psi\circ h^{-1})_*\nu_{h\cdot Z} = g_*\Psi_*\nu_Z$, and Theorem
\ref{naturality-thm} gives
${\cal B}(g_*\Psi_*\nu_Z) = g\cdot {\cal B}(\Psi_*\nu_Z) = g\cdot E\Psi(Z)$.

(c) $E(\operatorname{id})(Z) = {\cal B}(\nu_Z)=Z$ by Theorem
\ref{naturality-thm}, and the second statement follows from (b) with
$h=\operatorname{id}$, $\Psi = \operatorname{id}$.

(a) Consider
\[
{\cal G}(Z,W) = \int_{U(2)} \bigl(\Psi(U)-W\bigr)\bigl(1-W^*\Psi(U)\bigr)^{-1}
P(Z,U)\,dU ,
\]
which is real analytic on $\BB D^+\times\BB D^+$ (differentiation under the
integral sign is legitimate by the uniform bounds preceding Lemma
\ref{energy-cocycle-lem} together with \eqref{Poisson-Szego}) and which, by
Theorem \ref{barycenter-thm}, vanishes exactly when $W=E\Psi(Z)$. Fix
$Z_0 \in \BB D^+$; by (b) we may compose with elements of $U(2,2)$ and assume
$E\Psi(Z_0)=0$. Then, since $\frac{d}{ds}\bigl[(\Psi-sw)(1-sw^*\Psi)^{-1}
\bigr]_{s=0} = -w+\Psi w^*\Psi$,
\[
D_W{\cal G}(w)\big|_{(Z_0,0)} = -\,{\cal K}_{Z_0}(w), \qquad
{\cal K}_{Z_0}(w) := w - \int_{U(2)} \Psi(U)\,w^*\,\Psi(U)\,d\nu_{Z_0}(U) ,
\]
and, using $\tr\bigl((\Psi(U)w^*)^2\bigr)
= \overline{\tr\bigl((w\Psi(U)^*)^2\bigr)}$,
\[
\re\tr\bigl({\cal K}_{Z_0}(w)w^*\bigr)
= \|w\|_{HS}^2 - \re\tr\int_{U(2)}(wU'^*)^2\,d(\Psi_*\nu_{Z_0})(U')
= \tfrac12\,\frac{d^2}{dt^2}\Big|_{t=0}{\cal E}_{\Psi_*\nu_{Z_0}}(tw),
\]
which is strictly positive for $w \ne 0$ by \eqref{second-var} and Proposition
\ref{convexity-prop}. Hence $D_W{\cal G}$ is invertible at $(Z_0,0)$ and the
analytic implicit function theorem applies.
\end{proof}

\subsection{Consequences for the Schwarz Lemma}

We now show that the barycenter recovers a holomorphic self map of $\BB D^+$
that preserves the Shilov boundary from its boundary values, by a formula that
is natural for the full group $U(2,2)$ acting on both sides. We first record
an elementary mean value property; note that its proof uses the very device --
restriction to the holomorphic disks $t \mapsto tU_0$ -- on which the proof of
the Schwarz Lemma in Subsection \ref{Schwarz-subsect} is based.

\begin{lem}  \label{mean-value-lem}
Let ${\cal V}$ be a finite dimensional complex vector space and let
$f: \overline{\BB D^+} \to {\cal V}$ be continuous on $\overline{\BB D^+}$ and
holomorphic on $\BB D^+$. Then
\[
f(0) = \int_{U(2)} f(U)\,dU .
\]
\end{lem}

\begin{proof}
Fix $U_0 \in U(2)$. Since $\|tU_0\|=|t|$, the function $t \mapsto f(tU_0)$ is
holomorphic on ${\cal D}$ and continuous on $\overline{{\cal D}}$, so the
classical mean value property gives
$f(0) = \frac1{2\pi}\int_0^{2\pi}f(e^{i\vartheta}U_0)\,d\vartheta$.
Integrating this identity over $U_0 \in U(2)$ against $dU_0$ and using Fubini
together with the invariance of Haar measure under
$U \mapsto e^{i\vartheta}U$ gives the assertion.
\end{proof}

\begin{thm}  \label{F-is-extension-thm}
Let $F: \overline{\BB D^+} \to \overline{\BB D^+}$ be continuous on
$\overline{\BB D^+}$, holomorphic on $\BB D^+$, such that
$F(\BB D^+)\subseteq\BB D^+$ and $F\bigl(U(2)\bigr)\subseteq U(2)$.
Assume that $F_*\nu_Z$ is admissible for
every $Z \in \BB D^+$. Then
\begin{equation}  \label{F-bary}
F(Z) \;=\; {\cal B}\bigl(F_*\nu_Z\bigr) \qquad\text{for all }Z \in \BB D^+ ,
\end{equation}
that is, $F = E\bigl(F|_{U(2)}\bigr)$ is the barycentric extension of its own
restriction to the Shilov boundary. Equivalently, $W=F(Z)$ is the unique
solution in $\BB D^+$ of
\[
\int_{U(2)} \bigl(F(U)-W\bigr)\bigl(1-W^*F(U)\bigr)^{-1}\,P(Z,U)\,dU = 0 .
\]
\end{thm}

\begin{proof}
Fix $Z \in \BB D^+$ and put $P_1=F(Z) \in \BB D^+$. Choose $g,h \in U(2,2)$ with $g\cdot P_1=0$ and
$h\cdot 0 = Z$; by Lemma \ref{cocycle-lem} both act as real analytic
diffeomorphisms of $\overline{\BB D^+}$ preserving $U(2)$. Hence
\[
G := g\circ F\circ h : \overline{\BB D^+} \to \overline{\BB D^+}
\]
is continuous, holomorphic on $\BB D^+$, maps $U(2)$ into $U(2)$, and
satisfies $G(0)=g\cdot F(Z) = 0$. By Lemma \ref{mean-value-lem},
\[
0 = G(0) = \int_{U(2)}G(U)\,dU = \int_{U(2)}U'\,d\bigl(G_*\nu_0\bigr)(U') ,
\]
so $M(G_*\nu_0)=0$ and therefore ${\cal B}(G_*\nu_0)=0$ by Theorem
\ref{barycenter-thm}. Now $G_*\nu_0 = g_*F_*h_*\nu_0 = g_*F_*\nu_Z$ by
\eqref{nu-equivariance}, and $g_*F_*\nu_Z$ is admissible by Lemma
\ref{admissible-lem}(b). Applying Theorem \ref{naturality-thm},
\[
0 = {\cal B}\bigl(g_*F_*\nu_Z\bigr) = g\cdot{\cal B}\bigl(F_*\nu_Z\bigr),
\]
whence ${\cal B}(F_*\nu_Z) = g^{-1}\cdot 0 = P_1 = F(Z)$.
\end{proof}

\begin{cor}  \label{Schwarz-boundary-cor}
Let $F$ be as in Theorem \ref{F-is-extension-thm}. Then the following are
equivalent:
\begin{enumerate}
\item[(i)] $F(0)=0$;
\item[(ii)] $\displaystyle\int_{U(2)}F(U)\,dU = 0$;
\item[(iii)] ${\cal B}\bigl(F_*\nu_0\bigr)=0$.
\end{enumerate}
If they hold then, by Corollary \ref{Schwarz-cor0},
\[
\|F(Z)\| \le \|Z\| \quad \text{for all } Z \in \BB D^+,
\qquad \biggl\| \frac{\partial F}{\partial z_{ij}}(0)\biggr\| \le 1 .
\]
\end{cor}

\begin{proof}
(i)$\Leftrightarrow$(ii) is Lemma \ref{mean-value-lem}, and
(ii)$\Leftrightarrow$(iii) is the last statement of Theorem
\ref{barycenter-thm} applied to $\mu=F_*\nu_0$, since
$M(F_*\nu_0)=\int F(U)dU$.
\end{proof}

Thus the normalisation $F(0)=0$ required by the Schwarz Lemma can be verified
purely on the Shilov boundary. Here are two situations in which it holds
automatically.

\begin{cor}  \label{symmetry-cor}
Let $F$ be as in Theorem \ref{F-is-extension-thm} and suppose that the measure
$F_*\nu_0$ on $U(2)$ is invariant under $U \mapsto e^{i\vartheta_0}U$ for some
$\vartheta_0$ with $e^{i\vartheta_0}\ne1$ -- for instance, if $F|_{U(2)}$
preserves Haar measure, or if $F(e^{i\vartheta}U) = e^{ik\vartheta}F(U)$ for
some integer $k \ne 0$ and all $\vartheta$. Then $F(0)=0$, and consequently
$\|F(Z)\| \le \|Z\|$ on $\BB D^+$.
\end{cor}

\begin{proof}
Invariance gives $M(F_*\nu_0)=e^{i\vartheta_0}M(F_*\nu_0)$, hence
$M(F_*\nu_0)=0$; apply Corollary \ref{Schwarz-boundary-cor}. In the
homogeneity case, $\nu_0$ is invariant under $U\mapsto e^{i\vartheta}U$, so
$F_*\nu_0$ is invariant under $U \mapsto e^{ik\vartheta}U$ for every
$\vartheta$.
\end{proof}

\begin{ex}  \label{examples-ex}
The maps $F(Z)=aZ^kd$ with $a,d \in U(2)$ and $k \ge 1$, and
$F(Z) = N(Z)^k\cdot a$ with $a \in U(2)$, satisfy the hypotheses of Theorem
\ref{F-is-extension-thm} and of Corollary \ref{symmetry-cor}: they are
holomorphic on $\HC$, they map $U(2)$ into $U(2)$, they map $\BB D^+$ into
$\BB D^+$ (because $\|aZ^kd\|\le\|Z\|^k$ and $|N(Z)|=\lambda_1\lambda_2 \le
\|Z\|^2$), and they are homogeneous of degrees $k$ and $2k$ respectively.
Admissibility of $F_*\nu_Z$ follows from Lemma \ref{analytic-admissible-lem}.
Indeed, for $F(Z)=aZ^kd$ one has
$N(F(U)-W) = N(a)N(U^k - a^{-1}Wd^{-1})N(d)$, and choosing $U=\zeta\cdot1$
with $\zeta^k$ not an eigenvalue of $a^{-1}Wd^{-1}$ makes this nonzero. For
$F(Z)=N(Z)^ka$ one has, with $\zeta = N(U)$ on $U(2)$,
\[
N\bigl(F(U)-W\bigr) = \zeta^{2k}N(a) - \zeta^k \tr\bigl(a\operatorname{adj}W
\bigr) + N(W),
\]
a nonzero polynomial in $\zeta^k$ (its leading coefficient is $N(a)\ne0$),
which therefore vanishes for at most two values of $\zeta^k$.
\end{ex}

\begin{cor}  \label{rigidity-cor}
Let $F$ be as in Theorem \ref{F-is-extension-thm}. If $F|_{U(2)} = g|_{U(2)}$
for some $g \in U(2,2)$, then $F = g$ on $\BB D^+$. More generally, two maps
as in Theorem \ref{F-is-extension-thm} with the same restriction to $U(2)$
coincide.
\end{cor}

\begin{proof}
Immediate from \eqref{F-bary} and Theorem \ref{extension-thm}(c).
\end{proof}

\begin{rem}  \label{DE-rem}
Formula \eqref{F-bary} should be compared with the Poisson representation
\[F(Z)=\int_{U(2)}F(U)P(Z,U)\,dU,
\]
which also recovers $F$ from its boundary values.
The two are different: the Poisson integral is the {\em Euclidean}
mean of $F_*\nu_Z$, whereas ${\cal B}(F_*\nu_Z)$ is its {\em conformal}
barycenter. They agree here because the Euclidean mean of an admissible
measure vanishes exactly when its barycenter is the origin, and because both
sides are natural for $U(2,2)$. The advantage of the right hand side of
\eqref{F-bary} is that it makes sense, and produces a real analytic map
$\BB D^+ \to \BB D^+$ commuting with $U(2,2)$ on both sides, for boundary maps
$\Psi$ that need not extend holomorphically at all (Theorem
\ref{extension-thm}).
\end{rem}

\subsection{The Derivative of the Barycentric Extension}

The following formula is the analogue, in the present setting, of the implicit
differentiation of the barycenter equation that lies at the heart of the
method of \cite{BCG}. We use on $\HC$ the real inner product
$\langle x,y\rangle = \re \tr (xy^*)$, whose norm is $\|\cdot\|_{HS}$.

\begin{thm}  \label{derivative-thm}
Let $\Psi$ be as in Definition \ref{extension-df}, and suppose $E\Psi(0)=0$
(which can always be arranged by composing with an element of $U(2,2)$, by
Theorem \ref{extension-thm}(b)). Define real linear maps
${\cal K},{\cal L}: \HC \to \HC$ by
\[
{\cal K}(w) = w - \int_{U(2)} \Psi(U)\,w^*\,\Psi(U)\,dU,
\qquad
{\cal L}(v) = 4\int_{U(2)} \Psi(U)\;\re\tr\bigl(vU^*\bigr)\,dU .
\]
Then ${\cal K}$ is symmetric and positive definite for
$\langle\cdot,\cdot\rangle$, and
\begin{equation}  \label{dE-formula}
d(E\Psi)_0 = {\cal K}^{-1}\circ{\cal L} .
\end{equation}
Moreover $\|{\cal L}(v)\|_{HS} \le 4\|v\|_{HS}$, so that
\[
\bigl\| d(E\Psi)_0 \bigr\| \;\le\; \frac{4}{\kappa}, \qquad
\kappa := \min_{\|v\|_{HS}=1}
\Bigl( \|v\|_{HS}^2 - \re\tr \int_{U(2)}(vU'^*)^2\,d(\Psi_*\nu_0)(U')\Bigr)>0 .
\]
\end{thm}

\begin{proof}
Consider
\[
{\cal G}(Z,W) = \int_{U(2)} \bigl(\Psi(U)-W\bigr)\bigl(1-W^*\Psi(U)\bigr)^{-1}
P(Z,U)\,dU ,
\]
which is real analytic on $\BB D^+\times\BB D^+$ and vanishes exactly when
$W=E\Psi(Z)$, by Theorem \ref{barycenter-thm} and Definition
\ref{extension-df}. By hypothesis ${\cal G}(0,0)=0$.

Differentiating in $W$ at $(0,0)$ in the direction $w$: since
\[
\frac{d}{ds}\bigl[(\Psi-sw)(1-sw^*\Psi)^{-1}\bigr]_{s=0} = -w+\Psi w^*\Psi
\]
and $P(0,U)=1$, we get $D_W{\cal G}(w) = -{\cal K}(w)$.

Differentiating in $Z$ at $(0,0)$ in the direction $v$: at $W=0$ the integrand
is $\Psi(U)P(Z,U)$, and from
\[
\log P(Z,U) = 2\log\det(1-ZZ^*)-4\re\log\det(1-ZU^*)
\]
we obtain
\[
\frac{d}{ds}P(sv,U)\big|_{s=0} = 4\,\re\tr(vU^*);
\]
hence $D_Z{\cal G}(v)={\cal L}(v)$.

Symmetry of ${\cal K}$ is clear from
\[
\langle {\cal K}w_1,w_2\rangle = \re\tr(w_1w_2^*)
- \re\int\tr\bigl(\Psi(U) w_1^*\Psi(U)w_2^*\bigr)dU.
\]
Taking $w_1=w_2=w$ and
using $\tr\bigl((\Psi(U)w^*)^2\bigr) = \overline{\tr\bigl((w\Psi(U)^*)^2\bigr)}$
we get
\[
\langle {\cal K}w,w\rangle
= \|w\|_{HS}^2 - \re\tr\int_{U(2)}(wU'^*)^2\,d(\Psi_*\nu_0)(U')
= \tfrac12\,\frac{d^2}{dt^2}\Big|_{t=0}{\cal E}_{\Psi_*\nu_0}(tw) ,
\]
which is $>0$ for $w \ne 0$ by \eqref{second-var} and Proposition
\ref{convexity-prop}, since $\Psi_*\nu_0$ is admissible and $E\Psi(0)=0$.
Thus ${\cal K}$ is invertible and the analytic implicit function theorem gives
$D_Z{\cal G} + D_W{\cal G}\circ d(E\Psi)_0 = 0$, i.e. \eqref{dE-formula}.

For the bound on ${\cal L}$, note $\|\Psi(U)\|_{HS}=\sqrt2$ for
$\Psi(U) \in U(2)$, so by Cauchy-Schwarz
\[
\|{\cal L}(v)\|_{HS} \le 4\sqrt2 \int_{U(2)}\bigl|\tr(vU^*)\bigr|\,dU
\le 4\sqrt2\Bigl(\int_{U(2)}\bigl|\tr(vU^*)\bigr|^2dU\Bigr)^{1/2}
= 4\sqrt2 \cdot \tfrac1{\sqrt2}\|v\|_{HS} ,
\]
where we used the Schur orthogonality relations
$\int_{U(2)}U_{kl}\overline{U_{ij}}\,dU = \tfrac12\delta_{ik}\delta_{jl}$ for
the (irreducible, two dimensional) defining representation of $U(2)$, see
\cite{BD}.
\end{proof}

\begin{rem}  \label{derivative-check-rem}
For $\Psi = \operatorname{id}$ the two operators are computed from the Schur
relations together with $\int_{U(2)}U_{ij}U_{kl}\,dU = 0$ (which holds because
Haar measure is invariant under $U \mapsto e^{i\vartheta}U$, producing a factor
$e^{2i\vartheta}$): one finds ${\cal K}=\operatorname{id}$ and
${\cal L}=\operatorname{id}$, so \eqref{dE-formula} returns
$d(E\operatorname{id})_0=\operatorname{id}$, consistent with Theorem
\ref{extension-thm}(c).
\end{rem}

\begin{rem}  \label{Un-bary-rem}
As in Remark \ref{Un-n-rem}, this section carries over to
$\BB D^+_{(n)}=U(n,n)/(U(n)\times U(n))$ with Shilov boundary $U(n)$: one
takes the same $B_U(Z) = \log\bigl(|\det(1-ZU^*)|^2/\det(1-ZZ^*)\bigr)$ and
$P(Z,U)=e^{-nB_U(Z)}=\det(1-ZZ^*)^n/|\det(1-ZU^*)|^{2n}$, and
${\cal Z}_W = \{U \in U(n);\: \det(U-W)=0\}$.
Lemmas \ref{cocycle-lem}--\ref{admissible-lem}, \ref{energy-cocycle-lem},
\ref{expansion-lem}, \ref{trace-ineq-lem} and \ref{flats-lem}, and
Propositions \ref{convexity-prop} and \ref{coercivity-prop} together with
Theorems \ref{barycenter-thm}--\ref{derivative-thm}, all go through verbatim,
with one exception: Lemma \ref{Sigma-lem} was proved only for $n=2$, and it is
what shows that condition \eqref{star} alone forces the second variation
\eqref{second-var} to be non-degenerate. For general $n$ one should therefore
{\em define} $\mu$ to be admissible when both
\[
\sup_{W \in U(n)} \mu({\cal Z}_W)<\tfrac12
\qquad\text{and}\qquad
\mu(\Sigma_v)<1 \ \text{ for every } v \ne 0 ,
\]
where $\Sigma_v = \{U \in U(n);\: vU^*=Uv^*\}$; the first condition is used in
Proposition \ref{coercivity-prop} (note that $c_1+\dots+c_n \ge
(c_1^2+\dots+c_n^2)^{1/2}$ for $c_j\ge0$, which is what \eqref{slope-bound}
requires), the second in Proposition \ref{convexity-prop}. Both hold with
room to spare for every measure that is absolutely continuous with respect to
Haar measure, since each $\Sigma_v$, $v\ne0$, is the zero set of the real
analytic function $U \mapsto \|vU^*-Uv^*\|_{HS}^2$, which is not identically
zero (take $U = \zeta\cdot 1$ with $\zeta$ chosen so that $\zeta v$ is not
Hermitian; such $\zeta$ exists because $v$ and $iv$ cannot both be Hermitian).
\end{rem}

The barycenter method also gives a global estimate for maps into compact
quotients, distinct from the boundary construction above. The following
is an application of the higher-rank Degree Theorem of
Connell--Farb \cite{CF,CFE} to the irreducible rank-two symmetric space
$\BB D^+$.

\begin{thm}[Connell--Farb degree bound and consequences]
\label{Connell-Farb-thm}
Let $\Gamma<PU(2,2)=U(2,2)/Z(U(2,2))$ be a torsion-free cocompact lattice and
$M=\Gamma\backslash \BB D^+$. For every closed oriented Riemannian
$8$-manifold $(N,h)$ and continuous map $f:N\to M$, there is a constant
$C$, depending only on the dimension and the relevant lower
Ricci-curvature bounds, such that
\[
|\deg f|\,\operatorname{Vol}(M)
\le C\,\operatorname{Vol}(N,h).
\]
In particular, every continuous self-map $f:M\to M$ satisfies
$|\deg f|\le 1$, and
\[
\operatorname{Minvol}(M)>0,
\qquad\text{where}\qquad
\operatorname{Minvol}(M):=
\inf\{\operatorname{Vol}(M,g): |K_g|\le 1\},
\]
and $K_g$ denotes sectional curvature.
\end{thm}

\begin{rem}  \label{Connell-Farb-rem}
The self-map conclusion follows by applying the degree estimate to
iterates of $f$. For the minimal-volume conclusion, apply it to the
identity map from $(M,g)$ to $M$ with a fixed target metric. The
condition $|K_g|\le 1$ gives the required uniform lower Ricci-curvature bound.
\end{rem}

\section{Further Questions}  \label{questions-sect}

In this section we pose further questions that naturally arise from
our version of Schwarz Lemma but have not been addressed in this paper.

\begin{itemize}
\item
Given Theorem \ref{sqrt2-thm} and Proposition \ref{sqrt2-sharp-prop},
can we characterise the holomorphic maps $F$ for which the Lipschitz constant $\sqrt2$ is attained?

\item
The second part of the classical Schwarz Lemma asserts that if
$f: {\cal D} \to {\cal D}$ is a holomorphic function with $f(0)=0$ and
such that $|f(z_0)|=|z_0|$ for some $z_0 \in {\cal D}$ or $|f'(0)|=1$,
then $f$ is of the form $f(z) = \lambda z$, where $|\lambda|=1$. 
What conditions would force a holomorphic function $F: \BB D^+ \to \BB D^+$
be of the form $F(Z)=aZd^{-1}$ for some $a, d \in U(2)$?

Slightly more generally, what conditions would force a holomorphic
function $F: \BB D^+ \to \BB D^+$ to be a fractional linear transformation
coming from an element of $U(2,2)$?

Theorem \ref{infinitesimal-rigidity-thm} gives one answer when $F(0)=0$ and
$dF_0$ preserves the infinitesimal Kobayashi norm in every direction.  It is
natural to ask how far that hypothesis can be weakened, for instance whether
it is enough to impose equality only on a suitable family of rank-one
directions. Corollary \ref{rigidity-cor} gives a different boundary
criterion: if $F$ extends
continuously to $\overline{\BB D^+}$, maps $U(2)$ into $U(2)$, and its
restriction to $U(2)$ agrees with that of an element of $U(2,2)$, then $F$
is that element.


\item
Theorem \ref{Burns-Krantz-thm} gives a Schwarz-type rigidity theorem at the
smooth part of $\partial\BB D^+$.  What is the corresponding sharp one-point
statement at the nonsmooth Shilov boundary $U(2)$?  In particular, what jet
condition at a point of $U(2)$ forces a holomorphic self-map of $\BB D^+$ to
be the identity?

\item
Is there an analogue of Schwarz--Pick--Ahlfors Lemma \cite{A} for $\BB D^+$?
  
\item
Are there any other consequences of the classical Schwarz Lemma that
can be extended to $\BB D^+$?
  
\item
Can this construction be generalized to higher dimensions,
possibly as a homogeneous space $U(n,n)/(U(n) \times U(n))$?
(See Remarks \ref{Un-n-rem} and \ref{Un-bary-rem}: the sharp Lipschitz
constant becomes $\sqrt n$, and the barycenter theory goes through verbatim.)

\item
Is the barycenter map ${\cal B}$ of Theorem \ref{barycenter-thm} continuous
for the weak topology on $\{\mu ;\: \theta(\mu)\le\theta_0\}$, for each fixed
$\theta_0<\frac12$? Equivalently, is the coercivity estimate of Proposition
\ref{coercivity-prop} uniform over such families?

\item
Does the barycentric extension $E\Psi$ of Definition \ref{extension-df}
extend continuously to $\overline{\BB D^+}$ with boundary value $\Psi$, when
$\Psi$ is a homeomorphism of $U(2)$? In the one-variable case this is the
main technical theorem of \cite{DE}. A first step is to decide whether
${\cal B}(\mu_n) \to W_0$ whenever $\mu_n \to \delta_{W_0}$ weakly.

\item
Is $E\Psi$ a diffeomorphism of $\BB D^+$ when $\Psi$ is a homeomorphism of
$U(2)$? By Theorem \ref{derivative-thm} this amounts to showing that
${\cal K}^{-1}{\cal L}$ is invertible at every point.

\item
Can one find the sharp bound for $\bigl|\det d(E\Psi)_Z\bigr|$ in terms of
the two quadratic forms of Theorem \ref{derivative-thm}, in the spirit of the
sharp Jacobian estimates of \cite{BCG}? Note that the rank of $\BB D^+$ is
two, so that the Hessian \eqref{second-var} degenerates along the sets
$\Sigma_v$ of Lemma \ref{Sigma-lem}; the extremal measures for such an
estimate should be the ones concentrated near those sets.

\item
The Busemann function \eqref{Busemann-eqn} is, up to the Cayley transform of
Section \ref{half-plane-sect}, minus the logarithm of the determinant of the
``imaginary part'' (Remark \ref{Poisson-remarks}(ii)). Is there a barycenter
theory on the tube domain $\BB T^+$ that interacts with the causal structure
of the Minkowski space $\BB M$ and the cone $\C^+$?

\item
Is there a version of the barycenter construction for {\em regular} rather
than merely holomorphic functions, based on the Cauchy-Fueter kernel instead
of the Cauchy-Szeg\H o kernel?

\item
What are the domains of regularity for quaternionic left (or right) regular functions?
If $F$ is a left (or right) regular function defined on an open subset $\Omega_{\BB R} \subset \BB H$, then, by Lemma \ref{reg-ext-lem}, $F$ has an extension to a holomorphic left (or right) regular function defined on an open subset $\Omega \subset \HC$ containing $\Omega_{\BB R}$.
Furthermore, if, for example, ${\cal B}^{\epsilon} \subset \HC$ is the $\epsilon$-neighborhood of 
${\cal B} = \{ X \in \BB H ;\: \|X\|_{\BB H}<1 \}$ with $0<\epsilon<1/2$, then every holomorphic left (or right) regular function on ${\cal B}^{\epsilon}$ extends uniquely to a holomorphic left (or right) regular function on ${\cal B}^{\epsilon} \cup \BB D^+$. Thus we observe a quaternionic analogue of the Hartogs's phenomenon.

By analogy with domains of holomorphy in complex analysis \cite{Kr1}, let us call a connected open set $\Omega \subset \HC$ a {\em domain of regularity} if there do not exist non-empty open sets $\Omega_1, \Omega_2 \subset \HC$, with $\Omega_2$ connected, $\Omega_2 \not\subset \Omega$ and $\Omega_1 \subset \Omega \cap \Omega_2$, such that for every holomorphic left regular function $F$ on $\Omega$ there exists a holomorphic left regular function $\tilde F$ on $\Omega_2$ with $F = \tilde F$ on $\Omega_1$.
It is interesting to see if there is a description of domains of regularity in $\HC$ similar to domains of holomorphy.
\end{itemize}

\separate

\separate

\noindent
{\em Department of Mathematics, Indiana University,
Rawles Hall, 831 East 3rd St, Bloomington, IN 47405}

\end{document}